\documentclass[german,a4paper,12pt]{article}
\usepackage{amsmath,amsthm,amsfonts,latexsym,amscd}
\usepackage{amsthm}
\usepackage{amssymb}
\swapnumbers
\theoremstyle{plain}
\newtheorem{theorem}{Theorem}[section]
\newtheorem{theoreme}[theorem]{Th\'eor\`eme}
\newtheorem{lemme}[theorem]{Lemme}

\newtheorem{corollaire}[theorem]{Corollaire}
\newtheorem{proposition}[theorem]{Proposition}

\theoremstyle{definition}

\theoremstyle{remark}

\enlargethispage{-1cm}
\begin{document}

{\bf Extension de fibr\'es vectoriels et profondeur}\\

{\bf Helmut A. Hamm (M\"unster)}\\

\section{\bf Introduction}\par

Dans SGA2, A.Grothendieck a d\'emontr\'e des th\'eor\`emes de Lefschetz par voie alg\'ebrique, en particulier pour le groupe fondamental (alg\'ebrique) et le groupe de Picard [G2]. Les outils essentiels sont: conditions de profondeur, cohomologie locale, th\'eor\`emes de finitude et d'annulation,... Un des ingr\'edients est le th\'eor\`eme d'annulation suivant, en fixant un corps $k$:

\begin{theoreme} ([G2] XII Cor. 1.4) Soit $X$ un sch\'ema projectif sur $k$, $\cal S$ un faisceau alg\'ebrique coh\'erent sur $X$, ${\rm prof}\,{\cal S}\ge n$. Soit $\cal L$ un faisceau tr\`es ample sur $X$. Alors $H^q(X,{\cal S}\otimes_{{\cal O}_X} {\cal L}^{-l})=0$ pour $q<n$, $l\gg 0$.
\end{theoreme}

Ici, ${\rm prof}\,{\cal S}\ge n$ signifie que ${\rm prof}\,{\cal S}_x\ge n$ pour tout point ferm\'e $x$ de $X$, o\`u ${\rm prof}\,{\cal S}_x$ est la profondeur du ${\cal O}_{X,x}$-module ${\cal S}_x$. Naturellement il suffit de supposer que $\cal L$ est ample.

\vspace{2mm} 
Ce th\'eor\`eme \'etait utilis\'e, par example, afin de d\'emontrer un th\'eor\`eme de Lefschetz pour le groupe de Picard d'une vari\'et\'e projective ([G2] XII Cor. 3.6). En partie les r\'esultats y concernent, plus g\'en\'eralement, des fibr\'es vectoriels.

\vspace{2mm}
Le but de l'article pr\'esent est le passage du cas projectif au cas quasi-projectif - \`a part du traitement du cas analytique complexe correspondent. On ne va que p\'eparer ici le traitement du groupe de Picard et se concentrer sur le cas des fibr\'es vectoriels: une situation plus g\'en\'erale o\`u les r\'esultats sont plus fragmentaires.

\vspace{2mm}
En fait, on trouve d\'ej\`a un passage sur le cas quasi-projectif chez Grothendieck ( [G2], XII.5), mais en ce qui concerne notre genre de question il n'y a qu'une remarque qui se borne au cas o\`u la vari\'et\'e quasi-projective est le compl\'ementaire d'un nombre fini de points (loc.cit., 5.6).

\vspace{2mm}
Le passage du cas d'une vari\'et\'e projective \`a celui d'un sch\'ema quasi-projectif n\'ecessite  une condition de finitude. On va choisir une formulation de celle-ci qui \'etait utilis\'ee dans le cas complexe analytique par Trautmann et Siu, voir e.g. [BS]. Ceci marche dans notre contexte parce que nous avons un sch\'ema de Jacobson, voir [G3], comment on va le voir. On d\'efinit donc, pour un faisceau alg\'ebrique coh\'erent $\cal S$ sur un sch\'ema $Y$:
$$S_m({\cal S})=\{x\in Y\,\hbox{ferm\'e}\,|\, {\rm prof}\,{\cal S}_x\leq m\}.$$
On constate qu'en fait le th\'eor\`eme en haut s'ins\`ere dans un th\'eor\`eme de finitude g\'en\'eral, en passant de la vari\'et\'e projective au c\^one \'epoint\'e correspondant. Ainsi le th\'eor\`eme d'annulation en haut se g\'en\'eralise comme suit:

\newpage

\begin{theoreme} Soit $k$ un corps, $X$ un sch\'ema projectif sur $k$, $Y$ un sous-espace alg\'ebrique ferm\'e, $\cal S$ un faisceau alg\'ebrique coh\'erent sur $X\setminus Y$, ${\cal L}$ un faisceau ample sur $X$. Supposons que ${\rm prof}\,{\cal S}\ge m$ et que $\cal S$ v\'erifie la condition suivante de finitude $F_m$:\\
$\dim Y\cap \overline{S_l({\cal S})}< l-m$ pour tout $l\le m+\dim Y$.\\
Alors $H^s(X\setminus Y, {\cal S}\otimes_{{\cal O}_{X\setminus Y}}{\cal L}^{-\nu}|X\setminus Y)=0$, $s<m, \nu\gg 0$.
\end{theoreme}

Dans cet article nous posons $\dim\emptyset:=-\infty$. Pour la condition $F_m$ voir les conditions \'equivalentes de Proposition 2.1.

\vspace{2mm}
Dans section 4 on va appliquer ce genre de r\'esultats au semi-anneau $Vect\,X$ des classes d'isomorphie de fibr\'es vectoriels alg\'ebriques sur $X$. En particulier on va obtenir le th\'eor\`eme suivant:

\begin{theoreme} Soit $X\subset\mathbb{P}_N(k)$ un sous-sch\'ema projectif sur $k$, $Y$ une partie ferm\'ee, $H$ un hyperplan d\'efini par un faisceau d'id\'eaux $\cal I$ tel que ${\cal I}\otimes{\cal O}_X\simeq {\cal IO}_X$. Soit $\hat{X}$ la compl\'etion de $X$ le long de $X\cap H$. Supposons $prof\,{\cal O}_{X\cap H\setminus Y}\ge 2$ et que $\dim \overline{S_l({\cal O}_{X\setminus Y})}\cap Y\cap H<l-3$ pour tout $l\le \dim Y\cap H+3$. Alors
$$\lim\limits_{\longrightarrow} Vect\,U\simeq Vect\,(\hat{X}\setminus\hat{Y})\simeq \lim\limits_{\longleftarrow} Vect(X_n\setminus Y_n)$$
o\`u $U$ parcourt les voisinages ouverts de $X\cap H\setminus Y$ dans $X\setminus Y$, $n\in\mathbb{N}$ et $X_n$ est le $n$-i\`eme voisinage infinit\'esimal de $X\cap H$ dans $X$.
\end{theoreme}

On peut donc approximer $Vect\,\hat{X}\setminus\hat{Y}$ de deux c\^ot\'es. 

\vspace{2mm}
Il y a un cas o\`u on a une information un peu plus pr\'ecise:

\begin{theoreme} Soit ${\rm codim}_XY\ge 3$ et $X\setminus Y$ une intersection compl\`ete dans $\mathbb{P}_N(k)\setminus Y$, c.-\`a d. peut \^etre d\'efinie par $N\setminus \dim X$ \'equations, $\dim\,X\setminus Y\ge 3$. Soit $\cal E$ un fibr\'e vectoriel sur $X\setminus Y$. Alors $\cal E$ est trivial si et seulement si la restriction alg\'ebrique \`a $X\cap H\setminus Y$ est trivial.
\end{theoreme} 

Ce th\'eor\`eme permet de restreindre la question de trivialit\'e au cas de fibr\'es vectoriels sur une surface.

\vspace{2mm}

Il y a des renseignements plus pr\'ecis pour le groupe de Picard. En particulier on va d\'eduire dans section 5 le th\'eor\`eme de type de Lefschetz suivant:

\begin{theoreme} Soit $X$ un sous-sch\'ema projectif de $\mathbb{P}_N(k)$, $Y$ un ferm\'e de Zariski dans $X$ de codimension $\ge 4$, $X\setminus Y$ une intersection compl\`ete dans $\mathbb{P}_N(k)\setminus Y$ de dimension $\ge 3$. Alors $Pic(X\setminus Y)\simeq \mathbb{Z}$.
\end{theoreme}

Dans le cas $k=\mathbb{C}$ nous pouvons comparer avec le cadre analytique. Soit $X^{an}$ l'espace analytique complexe qui correspond au sch\'ema projectif complexe $X$.

\begin{theoreme} Soit $X\subset\mathbb{P}_N(\mathbb{C})$ un sous-sch\'ema projectif, $Y$ une partie ferm\'ee, $H$ un hyperplan d\'efini par un faisceau d'id\'eaux $\cal I$ tel que ${\cal I}\otimes{\cal O}_X\simeq {\cal IO}_X$. Soit $\dim S_{l+2}({\cal O}_{X\cap H\setminus Y})\le l$ pour $l\le \dim Y\cap H$. Alors on a un diagramme commutatif:

\footnotesize
$$\begin{array}{ccccccc}
\lim\limits_\rightarrow Vect\,U&&\stackrel{\simeq}{\longrightarrow}&&Vect(\hat{X}\setminus\hat{Y})&\stackrel{\simeq}{\rightarrow}&\lim\limits_\leftarrow Vect(X_n\setminus Y_n)\\
\downarrow\simeq&&&&\downarrow\simeq&&\downarrow\simeq\\
\lim\limits_\rightarrow Vect\,U^{an}&\stackrel{\simeq}{\rightarrow}&H^1(X^{an}\cap H^{an}\setminus Y^{an},Gl({\cal O}_{X^{an}}))&\stackrel{\simeq}{\rightarrow}&Vect(\hat{X}^{an}\setminus\hat{Y}^{an})&\stackrel{\simeq}{\rightarrow}&\lim\limits_\leftarrow Vect(X_n^{an}\setminus Y_n^{an})
\end{array}$$

\normalsize
o\`u $U$ parcourt les voisinages (de Zariski) de $X\cap H\setminus Y$ dans $X\setminus Y$.
\end{theoreme}

\begin{theoreme} Soient $X,Y,H$ comme dans Th\'eor\`eme 1.6. Supposons ${\rm codim}_XY\ge 3$, ${\rm codim}_XSing(X\setminus Y)\ge 3$, $\dim\,S_l({\cal O}_{X\setminus Y})\le l-2$ pour $l<\dim\,X$. Soit $\cal E$ un fibr\'e vectoriel sur $X\setminus Y$. Alors $\cal E$ est trivial si et seulement si la restriction analytique de ${\cal E}^{an}$ \`a $X^{an}\cap H^{an}\setminus Y^{an}$ est triviale.
\end{theoreme}

Ce th\'eor\`eme permet pour la question de trivialit\'e une r\'eduction au cas d'une surface projective lisse.

\vspace{2mm}
La situation est meilleure quand on regarde des fibr\'es vectoriels \`a connexion int\'egrable. Soit $Vect_{ci}(X)$ l'ensemble des classes d'isomorphie de fibr\'es vectoriels holomorphes \`a connexion int\'egrable, voir [D]. Alors 

\begin{theoreme}
Soit $X\subset\mathbb{P}_N(\mathbb{C})$ un sous-sch\'ema projectif, $Y$ une partie ferm\'ee, $X\setminus Y$ lisse, partout de dimension 
$\ge 3$, $H$ un hyperplan transverse \`a $X\setminus Y$, ${\rm codim}_XY\cap H\ge 3$. Alors on a un diagramme commutatif d'isomorphies:
$$\begin{array}{ccc}
Vect_{ci}(X\setminus Y)&\simeq&Vect_{ci}(X\cap H\setminus Y)\\
\downarrow\simeq&&\downarrow\simeq\\
Vect_{ci}(X^{an}\setminus Y^{an})&\simeq&Vect_{ci}(X^{an}\cap H^{an}\setminus Y^{an})
\end{array}$$
Notons qu'ici une connexion int\'egrable sur un fibr\'e vectoriel sur $X\setminus Y$ ou bien $X\cap H\setminus Y$ est automatiquement r\'eguli\`ere.
\end{theoreme}

\vspace{2mm}
Le travail pr\'esent\'e ici est dans la partie purement alg\'ebrique bas\'ee sur une lecture d\'etaill\'ee du travail admirable de Grothendieck [G2].

\vspace{2mm}
La partie sur le cas analytique et la comparaison alg\'ebrique/analytique a des racines compl\`etement diff\'erentes: le travail de Trautmann et Siu. Il y a ici des techniques d'extension li\'ees \`a une figure de Hartogs, apr\`es tout on constate qu'elles sont appliquables. Des m\'ethodes semblables sont d\'ej\`a utilis\'ees afin d'obtenir des th\'eor\`emes du type de Zariski-Lefschetz [Ha1] et pour le groupe de Picard local [Ha2] ou les groupes de Chow [Ha3]. Certains r\'esultats de cet article sont utilis\'es pour le groupe de Picard des vari\'et\'es projectives [HL2]. Notons que dans le cas complexe on dispose d'autres m\'ethodes transcendantes encore.

\section{\bf Conditions de finitude et d'annulation}

{\bf a) Conditions de finitude}

\vspace{2mm}

La possibilit\'e de travailler avec nos conditions de profondeur est bas\'ee sur l'observation suivante:

\begin{lemme} Supposons que $X$ est un sch\'ema de Jacobson (voir [G3] I 6.4.1). Soit $\cal S$ un faisceau alg\'ebrique coh\'erent sur $X$, $x\in X$. Alors:\\
a) ${\cal S}_x\neq 0$ $\Leftrightarrow$ ${\cal S}_z\neq 0$ pour tous les points ferm\'es $z$ de $\overline{\{x\}}$,\\
b)  Si $X$ est localement sous-sch\'ema d'un sch\'ema r\'egulier: prof ${\cal S}_x\le l-\dim\,\overline{\{x\}}$ $\Leftrightarrow$ tous les points ferm\'es de $\overline{\{x\}}$ appartiennent \`a $S_l({\cal S})$ $\Leftrightarrow$ $prof\,{\cal S}_z\le l$ pour tout point ferm\'e $z$ de $\overline{\{x\}}$.
\end{lemme}

Rappelons que $S_l({\cal S})$ \'etait d\'efini dans l'introduction.\\

{\bf D\'emonstration:} a) Comme $\cal S$ est coh\'erent on sait que le support de $\cal S$ est ferm\'e (voir [G3] 0, 5.2.2). Donc ``$\Rightarrow$'' est \'evident.\\
``$\Leftarrow$'': $\overline{\{x\}}$ et $\overline{\{x\}}\cap supp\,{\cal S}$ sont des ferm\'es tels que l'intersection avec le sous-espace des points ferm\'es de $X$ coincide. Par [G3] 0, 2.8.1 on conclut que $\overline{\{x\}}=\overline{\{x\}}\cap supp\,{\cal S}$, donc ${\cal S}_x\neq 0$.\\
b) Il suffit de consid\'erer le cas o\`u $X$ est r\'egulier de dimension $n$. Soit $pd$ la dimension projective d'un anneau local. Alors:\\
$prof\,{\cal S}_x+pd\,{\cal S}_x=\dim\,Spec\,{\cal O}_{X,x}=n-\dim\,\overline{\{x\}}$, voir [H1] III Prop. 6.12 A.\\
Or, $pd\,{\cal S}_x=\min\{j\,|\,j\ge 0, (Ext^i_{{\cal O}_X}({\cal S},{\cal O}_X))_x=0$ pour $i>j\}$, voir [H1] III Ex.6.6. Notons que $pd\,{\cal S}_x\le n$.\\
Le faisceau $Ext^i_{{\cal O}_X}({\cal S},{\cal O}_X)$ est coh\'erent. Donc:\\
$z\in S_l({\cal S})$ pour tous les points ferm\'es $z$ de $\overline{\{x\}}$\\
$\Leftrightarrow prof\,{\cal S}_z\le l$ pour tous les points ferm\'es $z$ dans $\overline{\{x\}}$\\
$\Leftrightarrow pd\,{\cal S}_z\ge n-l$ pour tous les points ferm\'es $z$ dans $\overline{\{x\}}$\\
$\Leftrightarrow \oplus_{n-l\le i\le n} Ext^i_{{\cal O}_X}({\cal S},{\cal O}_X)_z\neq 0$ pour tous les points ferm\'es $z$ dans $\overline{\{x\}}$\\
$\Leftrightarrow \oplus_{n-l\le i\le n} Ext^i_{{\cal O}_X}({\cal S},{\cal O}_X)_x\neq 0$ (voir a))\\
$\Leftrightarrow pd\,{\cal S}_x\ge n-l \Leftrightarrow prof\,{\cal S}_x\le l-\dim\,\overline{\{x\}}$.\\

Soient $X$, $Y$ comme dans l'introduction et $j:X\setminus Y\to X$ l'inclusion.

\vspace{2mm}
Afin de d\'emontrer Th\'eor\`eme 1.2 rappelons des r\'esultats de finitude qui sont essentiellement dues \`a [T2], [G2]:

\begin{proposition} Soit $\cal S$ un faisceau alg\'ebrique coh\'erent sur $X\setminus Y$ et $j:X\setminus Y\longrightarrow X$ l'inclusion. Pour $m\ge 1$, les conditions suivantes sont \'equivalentes:\\
a) $\cal S$ satisfait \`a la condition $F_m$: $\dim Y\cap\overline{S_l({\cal S})}<l-m$ pour tout $l\le m+\dim Y$,\\
b) pour tout $x\in X\setminus Y$: ${\rm prof}\,{\cal S}_x>m-c(x)$ avec $c(x):={\rm codim}\{\overline{\{x\}}\cap Y,\overline{\{x\}})$,\\
c) pour $s<m$, $R^sj_*{\cal S}$ est coh\'erent.\\
Dans ce cas, $\dim H^s(X\setminus Y,{\cal S})<\infty$, $s<m$, pourvu que $X$ est propre sur $k$.
\end{proposition}

{\bf Remarque:} Il suffit de supposer que $X$ est un sch\'ema de type fini sur $k$: dans ce cas il s'agit d'un sch\'ema de Jacobson qui est localement sous-sch\'ema d'un sch\'ema r\'egulier. \\
On peut d'ailleurs supposer que $Y\subset\overline{supp\,{\cal S}}$; dans ce cas, $(F_m)$ implique que ${\rm codim}_XY> m$ (posons $l=\dim\,X$).\\

{\bf D\'emonstration:} $b) \Longleftrightarrow c)$: On peut supposer que ${\cal S}=\hat{\cal S}|X\setminus Y$ o\`u $\hat{\cal S}$ est coh\'erent sur $X$: d'abord, $j_*{\cal S}$ est quasi-coh\'erent, donc limite inductive de ses sous-faisceaux coh\'erents, voir [G3] I 6.9.15.  Notons que ${\cal H}^0_Y(\hat{\cal S})$ est coh\'erent: sous l'hypoth\`ese b) ceci d\'ecoule de [G2] VIII Cor. 2.3, sous l'hypoth\`ese c) on sait que $j_*{\cal S}=j_*j^*\hat{\cal S}$ est coh\'erent, donc ${\cal H}^0_Y(\hat{\cal S})$ aussi. Alors c) est \'equivalent \`a la condition que ${\cal H}^s_Y\hat{\cal S}$ est coh\'erent pour $s\le m$. On est donc ramen\'e \`a un th\'eor\`eme de finitude de Grothendieck: [G2] VIII Cor. 2.3.\\
$a) \Longrightarrow b)$: Soit $x\in X\setminus Y$, $\dim\overline{\{x\}}=s$. Supposons que ${\rm prof}\,{\cal S}_x\le m-c(x)$: Alors tous les points ferm\'es de $\overline{\{x\}}\setminus Y$ appartiennent \`a $S_{m-c(x)+s}({\cal S})$, donc ceux de $\overline{\{x\}}$ \`a $\overline{S_{m-c(x)+s}({\cal S})}$, par Lemme 2.1. Ainsi, $s-c(x)=\dim\overline{\{x\}}\cap Y\le \dim \overline{S_{m-c(x)+s}({\cal S})}\cap Y<s-c(x)$ par hypoth\`ese, contradiction.\\
$b)\Longrightarrow a)$: Soit $x\in X\setminus Y$ tel que les points ferm\'es de $\overline{\{x\}}$ appartiennent \`a $\overline{S_l({\cal S})}$. Il faut montrer que $\dim\overline{\{x\}}\cap Y< l-m$. Mais ${\rm prof}\,{\cal S}_x\le l-s$, o\`u $s:=\dim\overline{\{x\}}$, par Lemme 2.1. Donc $\dim\overline{\{x\}}\cap Y=s-c(x)=(m-c(x))+(s-m)<{\rm prof}\,{\cal S}_x+(s-m)$ par hypoth\`ese, donc $\dim\overline{\{x\}}\cap Y<(l-s)+(s-m)=l-m$.\\
La conclusion est une cons\'equence de c).

\vspace{2mm}

Supposons maintenant que $X$ est un sous-sch\'ema de $\mathbb{P}_r$. On a $X=Proj\,A$, $Y=Proj\,B$, o\`u $A,B$ sont des $k$-alg\`ebres gradu\'ees, on peut supposer $A_0=B_0=k$; soit $\tilde{X}:=Spec\,A$, $\tilde{Y}:=Spec\,B$, $\tilde{j}:\tilde{X}\setminus\tilde{Y}\longrightarrow\tilde{X}$ l'inclusion, $0$ le point qui correspond \`a l'id\'eal $\oplus_{k>0} A_k$. On a une projection $\pi:\tilde{X}\setminus\tilde{Y}\longrightarrow X\setminus Y$. Soit $\tilde{\cal S}:=\pi^*{\cal S}$.

\begin{theoreme} Soit $\cal L$ ample sur $X$ et $m>0$. Les conditions suivantes sont \'equivalentes:\\
a) ${\rm prof}\,{\cal S}\ge m$, et $\cal S$ satisfait $(F_m)$,\\
b) $\tilde{\cal S}$ satisfait $(F_m)$,\\
c) $R^s\tilde{j}_*\tilde{\cal S}$ est coh\'erent pour $s<m$,\\
d) $\oplus_{n\in\mathbb{Z}}H^s(X\setminus Y,{\cal S}(n))$ est un $k[Z_0,\ldots,Z_r]$-module de type fini, $s<m$,\\
e) pour tout $s<m$, on a:\\
($e_1$) $H^s(X\setminus Y,{\cal S}\otimes{\cal L}^n)=0$, $n\ll 0$,\\
($e_2$) $\dim H^s(X\setminus Y,{\cal S}\otimes{\cal L}^n)<\infty$, $n\in\mathbb{Z}$,\\
($e_3$) $\oplus_{n\ge 0}H^s(X\setminus Y,{\cal S}(n))$ est un $k[Z_0,\ldots,Z_r]$-module de type fini.
\end{theoreme}

{\bf D\'emonstration:} a) $\Leftrightarrow$ b): On a ${\rm prof}\,\tilde{\cal S}_x={\rm prof}\,{\cal S}_{\pi(x)}+1$. La condition $(F_m)$ pour ${\cal S}$ dit que $\dim Y\cap \overline{S_l({\cal S})}< l-m$ pour tout $l\le m+\dim Y$. Ceci signifie que $\dim \tilde{Y}\cap \overline{S_l(\tilde{\cal S})}\setminus\{0\}< l-m$ pour tout $l\le m+\dim\tilde{Y}$. La condition ${\rm prof}\,{\cal S}_x\ge m$ pour tout point ferm\'e de $X\setminus Y$ dit que $S_l(\tilde{\cal S})=\emptyset$, c.-\`a d. $\tilde{Y}\cap \overline{S_l(\tilde{\cal S})}=\emptyset$ pour $l\le m$.\\
b) $\Leftrightarrow$ c): voir Proposition 2.2.\\
c) $\Leftrightarrow$ d): On sait de toute fa\c{c}on que les faisceaux $R^s\tilde{j}_*\tilde{\cal S}$ sont quasi-coh\'erents. Comme $\tilde{X}$ est affine, on obtient donc que $H^0(\tilde{X},R^s \tilde{j}_*\tilde{\cal S})=H^s(\tilde{X}\setminus\tilde{Y},\tilde{\cal S})\simeq \oplus_{n\in\mathbb{Z}}H^s(X\setminus Y,{\cal S}(n))$. En g\'en\'eral, on sait que $H^0(\tilde{X},R^s \tilde{j}_*\tilde{\cal S})$ est un $H^0(\tilde{X},{\cal O}_{\tilde{X}})$-module de type fini si et seulement si $R^s\tilde{j}_*\tilde{\cal S}$ est coh\'erent, voir [H1] II Cor. 5.5. Or, $H^0(\tilde{X},{\cal O}_{\tilde{X}})=A$, et $A$ est quotient de $R:=k[Z_0,\ldots,Z_r]$. D'o\`u le r\'esultat.\\
d) $\Rightarrow$ e): Soit ${\cal L}^k$ tr\`es ample. L'\'equivalence a) $\Leftrightarrow$ d) montre qu'on a les m\^emes hypoth\`eses pour ${\cal S}\otimes{\cal L}^j$, $j=0,\ldots,k-1$, au lieu de ${\cal S}$. On peut donc supposer que ${\cal L}={\cal O}_X(1)$.\\
La condition $(e_1)$ est alors \'evidente. Comme $R$ est noeth\'erien, $M_n:=\oplus_{k\ge n}H^s(X\setminus Y,{\cal S}\otimes{\cal L}^k)$ est un $R$-module de type fini pour tout $n$, ce qui donne $(e_3)$ pour $n=0$. Le m\^eme vaut donc pour $M_n/M_{n+1}=H^s(X\setminus Y,{\cal S}\otimes{\cal L}^n)$ aussi. Or, il s'agit en m\^eme temps d'un module sur $R/R_1\simeq k$, donc d'un espace vectoriel sur $k$ de dimension finie, on obtient donc $(e_2)$.\\
e) $\Rightarrow$ d): clair.

\vspace{2mm}
{\bf D\'emonstration de Th\'eor\`eme 1.2:} voir l'implication $a)\Longrightarrow e_1)$ de Th\'eor\`eme 2.3.\\

En particulier, on obtient le corollaire de Th\'eor\`eme 1.2 suivant qui peut \^etre d\'emontr\'e beaucoup plus facilement par voie directe:

\begin{corollaire} Soit $X$ un sch\'ema projectif sur $k$, $Y$ un sous-espace alg\'ebrique ferm\'e, $\cal L$ un faisceau ample sur $X$, $\cal S$ un faisceau alg\'ebrique coh\'erent sur $X\setminus Y$, ${\rm prof}\,{\cal S}\ge n$. Alors $H^q(X\setminus Y,{\cal S}\otimes_{{\cal O}_{X\setminus Y}} ({\cal L}|X\setminus Y)^{-\nu})=0$ pour $q<n-\dim Y-1, \nu\gg 0$.
\end{corollaire}

En fait, il y a une extension coh\'erente de $\cal S$ \`a $X$, on peut donc appliquer la m\'ethode de la d\'emonstration de [HL2] Th\'eor\`eme 2.1.

\vspace{2mm}
Le lemme suivant d\'ecoule de Th\'eor\`eme 2.3 aussi, voir $a)\Longrightarrow d)$,  mais peut \^etre d\'emontr\'e directement, voir Appendice:

\begin{lemme} Soit $\dim \overline{S_l({\cal S})}\cap Y<l-m$ pour tout $l\le m+\dim Y$, ${\rm prof}\,{\cal S}\ge m$. Alors $\oplus_{n\in\mathbb{Z}} H^s(X\setminus Y,{\cal S}(n))$ est un $k[Z_0,\ldots,Z_r]$-module de type fini pour $s<m$.
\end{lemme}

\begin{lemme} Soit $\dim \overline{S_l({\cal S})}\cap Y<l-1$ pour tout $l\le 1+\dim Y$. Alors ${\cal S}\otimes {\cal L}^n$ est engendr\'e par $\Gamma(X\setminus Y,{\cal S}\otimes{\cal L}^n)$, $n\gg 0$, le dernier espace \'etant de dimension finie.
\end{lemme}

{\bf D\'emonstration:} D'apr\`es Proposition 2.2, $j_*{\cal S}$ est coh\'erent, donc $j_*{\cal S}\otimes {\cal L}^n$ est engendr\'e par $\Gamma(X,j_*{\cal S}\otimes {\cal L}^n)=\Gamma(X\setminus Y,{\cal S}\otimes {\cal L}^n)$, $n\gg 0$, parce que $\cal L$ est ample, l'espace \`a droite \'etant de dimension finie.\\

{\bf b) Conditions d'annulation}

\vspace{2mm}
On aura aussi besoin des conditions d'annulation qui sont essentiellement dues \`a Grothendieck:

\begin{proposition} Soit $\cal S$ un faisceau alg\'ebrique coh\'erent sur $X$. Les conditions suivantes sont \'equivalentes:\\
a) $prof\,{\cal S}_y\ge n$ pour tout $y\in Y$,\\
b) $\dim\,Y\cap S_{l+n}({\cal S})\le l$ pour tout $l$,\\
c) ${\cal H}_Y^i{\cal S}=0$ pour $i<n$.\\
Dans ce cas, $H^i_Y(X,{\cal S})=0$ pour $i<n$.
\end{proposition}

{\bf Remarque:} Il suffit encore de supposer que $X$ est localement un sch\'ema de type fini sur $k$. \\

{\bf D\'emonstration:} a) $\Leftrightarrow$ c): voir [G2] III Prop. 3.3.\\
a) $\Rightarrow$ b): Soit $y\in Y$ tel que les points ferm\'es de $\overline{\{y\}}$ appartiennent \`a $S_{l+n}({\cal S})$. Il suffit de d\'emontrer que $s:=\dim\,\overline{\{y\}}\le l$. Par Lemme 2.1, $prof\,{\cal S}_y\le l+n-s$. Donc $\dim\overline{\{y\}}=s\le l+n-prof{\cal S}_y\le l+n-n=l$.\\
b) $\Rightarrow$ a): Supposons $y\in Y$, $prof{\cal S}_y<n$: Soit $s:=\dim\overline{\{y\}}$. Par Lemme 2.1, les points ferm\'es de $\overline{\{y\}}$ appartiennent \`a $S_{n+s-1}({\cal S})\cap Y$. Donc $\dim\,S_{n+s-1}({\cal S})\cap Y\ge \dim\overline{\{y\}}=s$, en contradiction avec l'hypoth\`ese.\\
La cons\'equence est obtenue par une suite spectrale, voir [G2] I Th. 2.6.\\

La condition b) est due \`a Scheja [Sch] dans le cas analogue analytique.

\begin{lemme} Soit $j:X\setminus Y\to X$ l'inclusion, $\cal S$ un faisceau alg\'ebrique coh\'erent sur $X\setminus Y$, et supposons que $\cal S$ satisfait \`a $(F_1)$. Alors $j_*{\cal S}$ est coh\'erent, et $prof\,(j_*{\cal S})_y\ge 2$ pour $y\in Y$.
\end{lemme}

{\bf D\'emonstration:} D'apr\`es Proposition 2.2 ${\cal T}:=j_*{\cal S}$ est coh\'erent. On a ${\cal T}\simeq j_*j^*{\cal T}$, donc ${\cal H}^k_Y({\cal T})=0, k=0,1$, ce qui implique $prof\,{\cal T}_y\ge 2, y\in Y$, par Proposition 2.7.

\section{\bf Cons\'equences pour les faisceaux formels}

{\bf a) Compl\'etion de faisceaux}

\vspace{2mm}
Soit $k$ un corps, $X$ un sch\'ema projectif sur $k$, $Y$ un sous-espace alg\'ebrique ferm\'e, $H$ un hyperplan et $\cal I$ le faisceau d'id\'eaux qui correspond \`a $H$. Pour un faisceau alg\'ebrique $\cal S$ soit $\hat{\cal S}$ la compl\'etion de $\cal S$ le long de $X\cap H$, voir [G3] I.10. On peut d\'eduire du Th\'eor\`eme 1.2, voir [G2] XII Th\'eor\`eme 2.1 pour $Y=\emptyset$:

\begin{theoreme} Soit $\cal S$ un faisceau alg\'ebrique coh\'erent sur $X\setminus Y$. Supposons que la fl\`eche naturelle ${\cal S}\otimes{\cal I}\rightarrow {\cal I}{\cal S}$ est injective, ${\rm prof}\,{\cal S}\ge m$ et que $\cal S$ v\'erifie la condition $F_m$, c.-\`a.d. $\dim Y\cap \overline{S_l({\cal S})}< l-m$ pour tout $l\le m+\dim Y$. Alors $H^s(X\cap H\setminus Y, \hat{\cal S})\stackrel{\simeq}{\longrightarrow} \lim\limits_{\stackrel{\leftarrow}{n}} H^s(X_n\setminus Y_n,{\cal S}/{\cal I}^n{\cal S})$ pour $s\le m-1$, et $H^s(X\setminus Y, {\cal S})\longrightarrow H^s(X\cap H\setminus Y, \hat{\cal S})$ est bijectif, $s<m-1$, et injectif, $s=m-1$.
\end{theoreme}

{\bf D\'emonstration:} On a la suite exacte longue de cohomologie associ\'ee \`a
$$0\longrightarrow {\cal I}^\nu{\cal S}\longrightarrow {\cal S}\longrightarrow {\cal S}/{\cal I}^\nu{\cal S}\longrightarrow 0$$
D'apr\`es Th\'eor\`eme 1.2, $H^s(X\setminus Y, {\cal S})\longrightarrow H^s(X\cap H\setminus Y, {\cal S}/{\cal I}^\nu{\cal S})$ est bijectif, $s<m-1$, et injectif, $s=m-1$, $\nu\gg 0$. Pour $s\le m-1$, ceci implique $H^s(X\cap H\setminus Y, \hat{\cal S})\simeq\lim\limits_{\stackrel{\leftarrow}{\nu}}H^s(X\cap H\setminus Y, {\cal S}/{\cal I}^\nu){\cal S} $, par [G4] $0_{III}$ 13.3.1, voir aussi $0_{III}$ 13.1.2. Ceci donne la bijectivit\'e pour $s<m-1$; pour $s=m-1$, la composition des fl\`eches\\
 $H^s(X\setminus Y,{\cal S})\longrightarrow  H^s(X\cap H\setminus Y, \hat{\cal S})\longrightarrow H^s(X\cap H\setminus Y, {\cal S}/{\cal I}^\nu){\cal S}$\\
est injective, $\nu\gg 0$, donc la premi\`ere fl\`eche aussi.\\

{\bf b) Cohomologie des faisceaux formels}

\vspace{2mm}
Maintenant commen\c{c}ons par un faisceau coh\'erent $\cal S$ sur $\hat{X}\setminus\hat{Y}$.
Soit $X\subset \mathbb{P}_r$, $H$ un hyperplan, $\cal I$ le faisceau d'id\'eaux correspondant. On peut supposer que $H$ est d\'efini par $Z_0=0$. Supposons que la fl\`eche naturelle ${\cal S}\otimes{\cal I}\longrightarrow {\cal I}{\cal S}$ est injective. Soit $\cal L$ un faisceau ample sur $X$.

\begin{lemme} Soit $\dim\overline{S_l({\cal S}/{\cal I}{\cal S})}\cap Y<l-m$ pour tout $l\le m+\dim Y\cap H$, ${\rm prof}\,{\cal S}/{\cal I}{\cal S}\ge m$. Alors le syst\`eme $(\oplus_{n\in\mathbb{Z}}H^s(X\cap H\setminus Y,({\cal S}/{\cal I}^k{\cal S})(n))_k$ v\'erifie pour $s<m-1$ la condition de Mittag-Leffler, voir [G4] $0_{III}$ 13.1.2.
\end{lemme}

{\bf D\'emonstration:} Soit ${\cal S}_k:={\cal S}/{\cal I}^k{\cal S}$ et $s<m$. Alors
$\oplus_{n\in\mathbb{Z}}H^s(X\cap H\setminus Y,gr{\cal S}^\cdot(n))$ $\simeq
(\oplus_{n\in\mathbb{Z}}H^s(X\cap H\setminus Y,{\cal S}_0(n)))\otimes_{k[Z_1,\ldots,Z_r]}k[Z_0,\ldots,Z_r]$. D'apr\`es Lemme 2.5, $\oplus_{n\in\mathbb{Z}}H^s(X\cap H\setminus Y,{\cal S}_0(n))$ est de type fini sur $k[Z_1,\ldots,Z_r]$. Le module \`a gauche est donc fini sur $k[Z_0,\ldots,Z_r]$. On d\'eduit le r\'esultat en question en utilisant [G4]  $0_{III}$ 13.7.7, corrig\'e comme dans [G5] ${\rm Err}_{III}$ 24, p. 89; voir [G2] IX Th\'eor\`eme 2.1, XII Lemme 3.3.

\begin{lemme} Sous l'hypoth\`ese du lemme pr\'ec\'edant on a pour tout $k<m$:
$\dim H^k(\hat{X}\setminus\hat{Y},{\cal S})=\dim\,H^k(\hat{X}\setminus\hat{Y},{\cal S}/{\cal I}^n{\cal S})<\infty, n\gg 0$ (donc $\dim H^k(\hat{X}\setminus\hat{Y},{\cal S}\otimes{\cal L}^n)<\infty$, $n$ arbitraire).
\end{lemme}

{\bf D\'emonstration:} D'abord $\dim H^k(\hat{X}\setminus\hat{Y},{\cal S}/{\cal I}^n{\cal S})<\infty$, $k<m$, par Proposition 2.2. De plus: $H^k(\hat{X}\setminus\hat{Y},{\cal I}^n{\cal S}/{\cal I}^{n+1}{\cal S})\simeq$ $H^k(\hat{X}\setminus\hat{Y},({\cal S}/{\cal I}{\cal S})(-n))=0$, $k<m$, $n\gg 0$, par Th\'eor\`eme 1.2, donc $H^k(\hat{X}\setminus\hat{Y},{\cal S}/{\cal I}^{n+1}{\cal S})\longrightarrow H^k(\hat{X}\setminus\hat{Y},{\cal S}/{\cal I}^n{\cal S})$ est injectif et $\lim\limits_{\stackrel{\leftarrow}{\nu}} H^k(\hat{X}\setminus\hat{Y},{\cal S}/{\cal I}^\nu{\cal S})\longrightarrow H^k(\hat{X}\setminus\hat{Y},{\cal S}/{\cal I}^n{\cal S})$ aussi pour $k<m, n\gg 0$. La condition de Mittag-Leffler est donc v\'erifi\'ee pour $(H^k(\hat{X}\setminus\hat{Y},{\cal S}/{\cal I}^n{\cal S}))_{n\ge 0}$, et $H^k(\hat{X}\setminus\hat{Y},{\cal S})\simeq  H^k(\hat{X}\setminus\hat{Y},\lim\limits_{\stackrel{\leftarrow}{\nu}}{\cal S}/{\cal I}^\nu{\cal S})\simeq H^k(\hat{X}\setminus\hat{Y},{\cal S}/{\cal I}^n{\cal S})$, $n\gg 0$.

\begin{lemme} Supposons que $\dim \overline{S_l({\cal S}/{\cal I}{\cal S})}\cap Y< l-2$ pour tout $l\le \dim Y\cap H+2$ et que ${\rm prof}\,{\cal S}/{\cal IS}\ge 2$. Alors $\Gamma(\hat{X}\setminus\hat{Y},{\cal S}\otimes{\cal L}^k)\otimes _{\Gamma(\hat{X}\setminus\hat{Y},{\cal O}_{\hat{X}})}{\cal O}_{\hat{X}\setminus\hat{Y}}\longrightarrow {\cal S}\otimes{\cal L}^k$ est surjectif, $k\gg 0$.
\end{lemme}

{\bf D\'emonstration:} Il suffit de traiter le cas ${\cal L}={\cal O}_X(1)$. On s'aper\c{c}oit que l'hypoth\`ese de Lemme 3.2 est donn\'ee avec $m=2$. Soit $G_k:=\oplus_{n\in\mathbb{Z}}\Gamma(X\cap H\setminus Y,({\cal S}/{\cal I}^k{\cal S})(n))$. D'apr\`es Lemme 3.2, il y a un $k_0$ tel que pour tout $k\ge k_0$ on a $Im(G_k\longrightarrow G_1)=Im(G_{k_0}\longrightarrow G_1)$, donc les deux c\^ot\'es coincident avec $Im(\lim\limits_{\stackrel{\leftarrow}{k}} G_k\longrightarrow G_1)=$ $Im \oplus_{n\in\mathbb{Z}}(\Gamma(\hat{X}\setminus \hat{Y},{\cal S}(n))\longrightarrow \oplus_{n\in\mathbb{Z}}\Gamma(X\cap H\setminus Y,({\cal S}/{\cal I}{\cal S})(n)))$.\hfill(*)\\
Soit $n\gg 0$. D'apr\`es Lemme 2.6, $({\cal S}/{\cal I}^{k_0}{\cal S})(n)$ est engendr\'e par $\Gamma(X\cap H\setminus Y,({\cal S}/{\cal I}^{k_0}{\cal S})(n))$. Alors $({\cal S}/{\cal IS})(n)$ est engendr\'e par l'image de $\Gamma(X\cap H\setminus Y,({\cal S}/{\cal I}^{k_0}{\cal S})(n))$ dans $\Gamma(X\cap H\setminus Y,({\cal S}/{\cal IS})(n))$, donc par celle de $\Gamma(\hat{X}\setminus\hat{Y},{\cal S}(n))$ dans $\Gamma(X\cap H\setminus Y,({\cal S}/{\cal IS})(n))$, voir (*). Par cons\'equent, le module $({\cal S}/{\cal IS})(n)$ est engendr\'e par $\Gamma(\hat{X}\setminus\hat{Y},{\cal S}(n))$, donc ${\cal S}(n)$ aussi par le lemme de Nakayama.

\vspace{2mm}
On a un analogue de Th\'eor\`eme 1.2 pour les faisceaux formels:

\begin{proposition} Supposons que ${\rm prof}\,{\cal S}/{\cal I}{\cal S}\ge m$, $\dim\overline{S_l({\cal S}/{\cal I}{\cal S})}\cap Y<l-m$ pour $l\le m+\dim Y\cap H$. Alors\\
a) $H^k(\hat{X}\setminus\hat{Y},{\cal S}\otimes{\cal L}^n)=0$ pour $k<m, n\ll 0$,\\
b) $H^k(\hat{X}\setminus\hat{Y},{\cal S})\simeq H^k(X_n\setminus Y_n,{\cal S}/{\cal I}^n{\cal S})$ pour $k<m,n\gg 0$.
\end{proposition}

{\bf D\'emonstration:} Soit $k<m$.\\
a) On peut supposer ${\cal L}={\cal O}_X(1)$. On a ${\rm prof}\,{\cal S}/{\cal I}{\cal S}\ge m$ et que ${\cal S}/{\cal I}{\cal S}$ satisfait $F_m$. Soit $n\gg 0$. Par Th\'eor\`eme 1.2 nous savons que $H^k(\hat{X}\setminus\hat{Y},{\cal I}^n{\cal S}/{\cal I}^{n+1}{\cal S})=0$. Ceci implique que $H^k(\hat{X}\setminus\hat{Y},{\cal I}^n{\cal S}/{\cal I}^{n+\nu}{\cal S})=0$ pour $\nu\ge 0$. Le syst\`eme $(H^k(\hat{X}\setminus\hat{Y},{\cal I}^n{\cal S}/{\cal I}^{n+\nu}{\cal S}))$ v\'eriefie la condition de Mittag-Leffler, et $H^k(\hat{X}\setminus\hat{Y},{\cal I}^n{\cal S})=0$.\\
b) Gr\^ace \`a la d\'emonstration de Lemme 3.3, le syst\`eme $(H^k(\hat{X}\setminus\hat{Y},{\cal S}/{\cal I}^n{\cal S})$ v\'erifie la condition de Mittag-Leffler, et ${\cal S}\simeq \lim\limits_{\stackrel{\leftarrow}{n}}{\cal S}/{\cal I}^n{\cal S}$, donc
$H^k(\hat{X}\setminus\hat{Y},{\cal S})\simeq \lim\limits_{\stackrel{\leftarrow}{n}}H^k(\hat{X}\setminus\hat{Y},{\cal S}/{\cal I}^n{\cal S})$.\\
D'autre part, la fl\`eche $H^k(\hat{X}\setminus\hat{Y},{\cal S}/{\cal I}^{n+1}{\cal S})\longrightarrow H^k(\hat{X}\setminus\hat{Y},{\cal S}/{\cal I}^n{\cal S})$ est injective, $n\gg 0$, comme on a vu ci-dessus. Pour $n\gg 0$, on conclut que $H^k(\hat{X}\setminus\hat{Y},{\cal S}/{\cal I}^{n+1}{\cal S})\simeq H^k(\hat{X}\setminus\hat{Y},{\cal S}/{\cal I}^n{\cal S})$, d'o\`u le r\'esultat cherch\'e.

\vspace{2mm}
Notons que l'hypoth\`ese de Proposition 3.5 est v\'erifi\'ee si ${\cal S}=\hat{\cal G}$ et $\cal G$ v\'erifie l'hypoth\`ese de Th\'eor\`eme 3.1 avec $m+1$ au lieu de $m$.\\

{\bf c) Extension de faisceaux formels}

\vspace{2mm}
On a donc la g\'en\'eralisation suivante de [G2] XII Th\'eor\`eme 3.1 ($Y=\emptyset$):

\begin{theoreme} Sous l'hypoth\`ese de Lemme 3.4, il y a un faisceau coh\'erent $\cal G$ sur $X$ tel que $\hat{\cal G}|\hat{X}\setminus\hat{Y}\simeq {\cal S}$.
\end{theoreme}

{\bf D\'emonstration:} Nous pouvons supposer que $Y\subset \overline{supp\,{\cal S}}$ et $X=\mathbb{P}_r(k)$, alors $r\ge 2$ \`a cause de l'hypoth\`ese sur la profondeur. La remarque apr\`es Proposition 2.1 montre que ${\rm codim}_XY\ge 4$. Soit $n\gg 0$, alors ${\cal S}(n)$ est engendr\'e par $\Gamma(\hat{X}\setminus\hat{Y},{\cal S}(n))$ par Lemme 3.4, et le dernier espace est de dimension finie d'apr\`es Lemme 3.3. On a donc un \'epimorphisme ${\cal O}^d_{\hat{X}\setminus\hat{Y}}\longrightarrow {\cal S}(n)$, donc ${\cal O}^d_{\hat{X}\setminus\hat{Y}}(-n)\longrightarrow {\cal S}$. Soit $\cal K$ le noyau qui a les m\^emes propri\'et\'es que $\cal S$. On a donc un \'epimorphisme ${\cal O}^c_{\hat{X}\setminus\hat{Y}}\longrightarrow {\cal K}(l)$, donc ${\cal O}^c_{\hat{X}\setminus\hat{Y}}(-l)\longrightarrow {\cal K}$. La suite
${\cal O}^c_{\hat{X}\setminus\hat{Y}}(-l)\longrightarrow{\cal O}^d_{\hat{X}\setminus\hat{Y}}(-n)\longrightarrow {\cal S}\longrightarrow 0$ est donc exacte. L'homomorphisme ${\cal O}^c_{\hat{X}\setminus\hat{Y}}(-l)\longrightarrow{\cal O}^d_{\hat{X}\setminus\hat{Y}}(-n)$ correspond \`a un \'el\'ement de $\Gamma(\hat{X}\setminus\hat{Y},Hom({\cal O}^c_{\hat{X}\setminus\hat{Y}}(-l),{\cal O}^d_{\hat{X}\setminus\hat{Y}}(-n))\simeq \Gamma(\hat{X}\setminus\hat{Y},\hat{\cal T})$ avec\\
${\cal T}:=Hom({\cal O}^c_{X\setminus Y}(-l),{\cal O}^d_{X\setminus Y}(-n))$. Par Th\'eor\`eme 3.1, on obtient un \'el\'ement correspondant de $\Gamma(X\setminus Y,{\cal T})$. Ceci signifie que l'homomorphisme provient d'un homomorphisme ${\cal O}^c_{X\setminus Y}(-l)\longrightarrow{\cal O}^d_{X\setminus Y}(-n)$. Posons ${\cal G}_1:={\rm coker}({\cal O}^c_{X\setminus Y}(-l)\longrightarrow{\cal O}^d_{X\setminus Y}(-n))$ et $\cal G$ une extension coh\'erente de ${\cal G}_1$ \`a $X$. Alors $\hat{\cal G}|\hat{X}\setminus\hat{Y}\simeq {\cal S}$.

\begin{corollaire} Sous l'hypoth\`ese de Lemme 3.4 on peut \'etendre $\cal S$ \`a un faisceau coh\'erent sur $\hat{X}$.
\end{corollaire}

On peut essayer de d\'emontrer Th\'eor\`eme 3.6 comme suit: \'etendre d'abord \`a $\hat{X}$ et apr\`es \`a $X$. En fait, l'extension \`a $\hat{X}$ est possible, mais la deuxi\`eme \'etape est probl\'ematique, voir Appendice.

\section{\bf Cons\'equences pour les fibr\'es vectoriels alg\'ebriques}

Soit d'abord $X$ un sch\'ema de type fini sur $k$, $Y$ un sous-espace ferm\'e.

\begin{theoreme} Soient $Y$ un sous-espace ferm\'e de $X$, $\dim Y\cap S_{l+2}({\cal O}_X)\le l$ pour tout $l$. Alors $Vect\,X\longrightarrow Vect\,(X\setminus Y)$ est injectif.
\end{theoreme}

{\bf D\'emonstration:} Soient ${\cal E}$ et ${\cal E}'$ des faisceaux coh\'erents localement libres sur $X$ tels que ${\cal E}|X\setminus Y\simeq {\cal E}'|X\setminus Y$. On raisonne comme dans la d\'emonstration de Cor. 2.4 dans [G2] XII: Soit ${\cal G}:=Hom({\cal E},{\cal E}')$. L'isomorphisme donne une section de ${\cal G}|X\setminus Y$ qui vient d'une section de $\cal G$: l'hypoth\`ese garantit que la restriction $\Gamma(X,{\cal G})\longrightarrow \Gamma(X\setminus Y,{\cal G})$ est un isomorphisme, voir Proposition 2.7: on a $H^k_Y(X,{\cal G})=0, k=0,1$. La section de $\cal G$ d\'efinit un homomorphisme ${\cal E}\longrightarrow {\cal E}'$, on v\'erifie qu'il s'agit d'un isomorphisme: on prolonge l'isomorphisme inverse et conclut par unicit\'e qu'il s'agit de l'homomorphisme inverse \`a l'extension d\'ej\`a trouv\'ee.\\

Soit maintenant $X\subset \mathbb{P}_m(k)$ un sch\'ema projectif sur $k$, $Y$ un sous-espace alg\'ebrique ferm\'e, $H$ un hyperplan tel que ${\cal I}\otimes{\cal O}_X\simeq {\cal I}{\cal O}_X$, $\cal I$ \'etant le faisceau d'id\'eaux qui correspond \`a $H$. Soit $\hat{X}$ la compl\'etion de $X$ le long de $X\cap H$.

\begin{theoreme} Supposons $prof\,{\cal O}_{X\setminus Y}\ge 2$ et que ${\cal O}_{X\setminus Y}$ v\'erifie la condition $(F_2)$. Alors:\\
a) Si $\cal E$ est un faisceau coh\'erent localement libre sur $X\setminus Y$, on a des isomorphismes
$\Gamma(X\setminus Y,{\cal E})\simeq\Gamma(\hat{X}\setminus \hat{Y},\hat{\cal E})\simeq \lim\limits_{\stackrel{\leftarrow}{n}}\Gamma(X_n\setminus Y_n,{\cal E}/{\cal I}^n{\cal E})$.\\
b) $Vect\,(X\setminus Y)\longrightarrow \lim\limits_{\stackrel{\leftarrow}{n}} Vect\,(X_n\setminus Y_n)$ est injectif, donc  $Vect\,(X\setminus Y)\longrightarrow Vect\,(\hat{X}\setminus \hat{Y})$ aussi.
\end{theoreme}

{\bf D\'emonstration:} a) Ceci d\'ecoule de Th\'eor\`eme 3.1.\\
b) On raisonne comme dans la d\'emonstration de Th\'eor\`eme 4.1. Soient ${\cal F}_1$ et ${\cal F}_2$ des faisceaux coh\'erents localement libres sur $X\setminus Y$ tel que ${\cal F}_1|X_n\setminus Y_n\simeq {\cal F}_2|X_n\setminus Y_n$ pour tout $n$. Soit ${\cal G}:=Hom({\cal F}_1,{\cal F}_2)$. Par Th\'eor\`eme 1.2 il y a un $n$ tel que $H^0(X\setminus Y,{\cal G})\simeq H^0(X_n\setminus Y_n,{\cal G}|X_n\setminus Y_n)$. L'isomorphisme ${\cal F}_1|X_n\setminus Y_n\simeq {\cal F}_2|X_n\setminus Y_n$ donne une section de ${\cal G}|X_n\setminus Y_n$ qui vient d'une section de $\cal G$. Celle-ci d\'efinit un homomorphisme ${\cal F}_1\longrightarrow {\cal F}_2$, on v\'erifie qu'il s'agit d'un isomorphisme, comme dans la d\'emonstration de Th\'eor\`eme 4.1.\\

Rappelons les conditions $Lef$ et $Leff$ de Lefschetz introduites dans SGA2 ([G2] X 2, p. 112):

\vspace{2mm}
Soit $X$ un sch\'ema de type fini sur $k$, $Z$ une partie ferm\'ee de $X$. Pour un faisceau alg\'ebrique coh\'erent $\cal F$ sur $X$ soit $\hat{\cal F}$ sa compl\'etion le long de $Z$. Soit $\hat{X}$ la compl\'etion de $X$ le long de $Z$.

\vspace{2mm}
{\bf D\'efinition:} On a $Lef(X,Z)$ si pour tout ouvert $U$ dans $X$ tel que $Z\subset U$ et pour tout faisceau alg\'ebrique coh\'erent localement libre $\cal E$ sur $U$ on a $H^0(U,{\cal E})\simeq H^0(\hat{X},\hat{\cal E})$.\\
On a $Leff(X,Z)$ si l'on a $Lef(X,Z)$ et pour tout faisceau alg\'ebrique coh\'erent localement libre ${\cal E}'$ sur $\hat{X}$ il y a un voisinage ouvert $U$ de $Z$ dans $X$ et un faisceau coh\'erent localement libre $\cal E$ sur $U$ tels que $\hat{\cal E}\simeq {\cal E}'$.

\vspace{2mm}
Retournons au contexte habituel.

\begin{corollaire} Supposons $prof\,{\cal O}_{X\setminus Y}\ge 2$, $\dim \overline{S_l({\cal O}_{X\setminus Y})}\cap Y\cap H<l-3$ pour tout $l\le \dim Y\cap H+3$.\\
a) On a $Lef(X\setminus Y,X\cap H\setminus Y)$.\\
b) Pour tout voisinage ouvert $U$ de $X\cap H\setminus Y$ dans $X\setminus Y$, la fl\`eche $Vect\,U\to Vect\,\hat{X}\setminus\hat{Y}$ est injective.
\end{corollaire}

{\bf D\'emonstration:} a) Soit $U$ un voisinage ouvert de $X\cap H\setminus Y$ dans $X\setminus Y$. Posons $Y':=X\setminus U$. On s'aper\c{c}oit que les hypoth\`eses de Th\'eor\`eme 3.1 sont remplies avec $m=2$ et $Y'$, ${\cal O}_U$ au lieu de $Y$ et $\cal S$: on a $Y'\cap H=Y\cap H$. Notons que $S_l({\cal O}_{X\setminus Y'})\subset S_l({\cal O}_{X\setminus Y})$, donc $\dim \overline{S_l({\cal O}_{X\setminus Y'})}\cap Y'\le \dim\overline{S_l({\cal O}_{X\setminus Y})}\cap Y\cap H+1<l-2$, et $prof\,{\cal O}_{X\setminus Y'}\ge 2$.\\
b) par Th\'eor\`eme 4.2, avec $Y':=X\setminus U$ au lieu de $Y$.\\

La comparaison de $Vect\,(\hat{X}\setminus\hat{Y})$ et $Vect\,(X\cap H\setminus Y)$ est plus difficile. Notons qu'il ne s'agit pas des groupes mais des ensembles avec un \'el\'ement distingu\'e. On peut donc parler du noyau mais il faut \^etre prudent parce que la trivialit\'e du noyau n'implique pas l'injectivit\'e.
Le sous-ensemble de $Vect\ldots$ repr\'esent\'e par des fibr\'es vectoriels de rang $r$ sera d\'enot\'e par $Vect_r\ldots$. D'abord:

\begin{lemme}
Soit $n\in\mathbb{N}$ tel que $H^1(X\cap H\setminus Y,{\cal I}^k/{\cal I}^{k+1})=0$ pour tout $k\ge n\ge 1$. Alors le noyau de $\lim\limits_{\stackrel{\leftarrow}{m}}Vect\,(X_m\setminus Y_m)\longrightarrow Vect\,(X_n\setminus Y_n)$ est trivial.
\end{lemme}

{\bf D\'emonstration:} Soit $r\in\mathbb{N}$. De fa\c{c}on analogue \`a [G2] XI (1.1), on a une suite exacte:
$$0\longrightarrow ({\cal I}^k/{\cal I}^{k+1})^{\oplus r^2}\longrightarrow GL_r({\cal O}_{X_{k+1}})\longrightarrow GL_r({\cal O}_{X_k})\longrightarrow 1$$
On passe \`a la suite exacte de cohomologie, cf. [Fr], qui se termine ici par $H^1(X\cap H\setminus Y,GL_r({\cal O}_{X_k}))$, et d\'eduit que le noyau de $Vect_r\,(X_{k+1}\setminus Y_{k+1})\longrightarrow Vect_r\,(X_k\setminus Y_k)$ est trivial, $k\ge n$, donc celui de $\lim\limits_{\stackrel{\leftarrow}{m}}Vect_r\,(X_m\setminus Y_m)\longrightarrow Vect_r\,(X_n\setminus Y_n)$ aussi.

\vspace{2mm}
En particulier, on obtient en utilisant Th\'eor\`eme 1.2, parce que ${\cal I}^k/{\cal I}^{k+1}\simeq ({\cal I}/{\cal I}^2)^k$:

\begin{proposition} Supposons que $prof\,{\cal O}_{X\cap H\setminus Y}\ge 1$ et que ${\cal O}_{X\cap H\setminus Y}$ v\'erifie la condition $(F_1)$. Pour $n\gg 0$. le noyau de  $\lim\limits_{\stackrel{\leftarrow}{m}}Vect\,(X_m\setminus Y_m)\longrightarrow Vect\,(X_n\setminus Y_n)$ est trivial.
\end{proposition}

Sous des hypoth\`eses plus fortes on peut appliquer Lemme 4.4 avec $n=1$:

\begin{theoreme} Soit ${\rm codim}_XY\cap H\ge 4$ et $X\setminus Y$ une intersection compl\`ete dans $\mathbb{P}_N(k)\setminus Y$ de dimension $\ge 3$. Alors le noyau de $\lim\limits_{\stackrel{\leftarrow}{m}}Vect\,(X_m\setminus Y_m)\longrightarrow Vect\,(X\cap H\setminus Y)$ est trivial.
\end{theoreme}

{\bf D\'emonstration:} Soit $n\ge 1$. Il suffit de d\'emontrer que le noyau de $Vect(X_{n+1}\setminus Y_{n+1})\to Vect(X_n\setminus Y_n)$ est trivial, o\`u bien $H^1(X\setminus Y,{\cal I}^n/{\cal I}^{n+1})=0$. Soit $i:X\to\mathbb{P}_N$ l'inclusion. Il y a une r\'esolution de $(i_*({\cal I}^n/{\cal I}^{n+1}))|\mathbb{P}_N\setminus Y\cap H$ par des faisceaux qui sont sommes directes de ${\cal O}_{{\mathbb P}_N\setminus Y\cap H}(l)$ avec $l<0$, le nombre des faisceaux non-triviaux \'etant $\le k+2$, o\`u $k:=N-\dim\,X$. Il suffit donc de v\'erifier que $H^j(\mathbb{P}_N\setminus Y\cap H,{\cal O}_{\mathbb{P}_N}(l))=0$, $j\le k+2$, ce qui provient du fait que ${\cal H}^j_{Y\cap H}({\cal O}_{\mathbb{P}_N})=0$, donc $H^j_{Y\cap H}(\mathbb{P}_N,{\cal O}_{\mathbb{P}_N}(l))=0$, $j\le k+3$, voir Proposition 2.7, et $H^j(\mathbb{P}_N,{\cal O}_{\mathbb{P}_N}(l))=0$, $l<0, j\le k+2$.

\vspace{2mm}
{\bf D\'emonstration de Th\'eor\`eme 1.4:} d\'ecoule de Th\'eor\`eme 4.2 et 4.6.

\begin{theoreme} Supposons que $prof{\cal O}_{X\cap H\setminus Y}\ge 1$ et que ${\cal O}_{X\cap H\setminus Y}$ v\'erifie la condition $(F_1)$. Alors $Vect\,(\hat{X}\setminus \hat{Y})\simeq \lim\limits_{\stackrel{\leftarrow}{m}}Vect(X_m\setminus Y_m)$.
\end{theoreme}

{\bf D\'emonstration:} D'abord, on a la surjectivit\'e par [G3] I 10.10.8.6:\\
Soit $({\cal E}_n)$ un repr\'esentant d'un \'el\'ement de $\lim\limits_{\longleftarrow} Vect\,(X_n\setminus Y_n)$, c.-\`a d. ${\cal E}_n$ est coh\'erent et localement libre sur $X_n\setminus Y_n$, ${\cal E}_{n+1}|(X_n\setminus Y_n)\simeq {\cal E}_n$ pour tout $n$. Alors il y a un faisceau $\cal E$ sur $\hat{X}\setminus\hat{Y}$ qui est coh\'erent et localement libre tel que ${\cal E}|X_n\setminus Y_n\simeq {\cal E}_n$ pour tout $n$.\\
Pour l'injectivit\'e, supposons que $\cal E$ et ${\cal E}'$ sont localement libres de type fini sur $\hat{X}\setminus\hat{Y}$ tels que ${\cal E}|X_n\setminus Y_n\simeq {\cal E}'|X_n\setminus Y_n$ pour tout $n$. Le probl\`eme est que les isomorphismes a priori ne doivent pas \^etre compatibles.\\
Posons ${\cal S}:=Hom({\cal E},{\cal E}')$. Ce faisceau est coh\'erent et localement libre. Choisissons $n$ suffisamment grand tel qu'on a d'abord $H^0(\hat{X}\setminus\hat{Y},{\cal S})\simeq H^0(\hat{X}\setminus\hat{Y},{\cal S}/{\cal I}^n{\cal S})$, par Proposition 3.5. On y pose $m=1$, ${\cal S}:={\cal O}_{\hat{X}}$.\\
On veut d'abord montrer que l'isomorphisme ${\cal E}|X_n\setminus Y_n\longrightarrow {\cal E}'|X_n\setminus Y_n$ s'\'etend \`a un morphisme ${\cal E}\longrightarrow {\cal E}'$: Le premier d\'efinit un \'el\'ement de $H^0(\hat{X}\setminus\hat{Y},{\cal S}/{\cal I}^n{\cal S})$, et ce groupe est isomorphe \`a $H^0(\hat{X}\setminus\hat{Y},{\cal S})$ par le choix de $n$.\\
Comme dans la d\'emonstration de Th\'eor\`eme 4.1 on montre qu'il s'agit d'un isomorphisme: on utilise  $Hom({\cal E}',{\cal E})$, $Hom({\cal E},{\cal E})$ et $Hom({\cal E}',{\cal E}')$ au lieu de ${\cal S}:=Hom({\cal E},{\cal E}')$.

\vspace{2mm}
{\bf D\'emonstration de Th\'eor\`eme 1.3:} L'isomorphisme \`a droite est assur\'e par le th\'eor\`eme pr\'ec\'edant. \\
Il reste \`a d\'emontrer la bijectivit\'e de la fl\`eche \`a gauche.\\
L'injectivit\'e d\'ecoule de Corollaire 4.3.\\
Surjectivit\'e: Soit $\cal E$ coh\'erent localement libre sur $\hat{X}\setminus \hat{Y}$, alors il y a un faisceau coh\'erent $\cal F$ sur $X\setminus Y$ tel que $\hat{\cal F}\simeq {\cal E}$, par Th\'eor\`eme 3.6. Comme $\hat{\cal F}$ est localement libre on conclut que ${\cal F}$ est localement libre sur un voisinage $U$, par [G3] I 10.8.15. Alors ${\cal F}|U$ repr\'esente l'image inverse cherch\'ee. Voir [G2] XII Corollaire 3.4b dans le cas $Y=\emptyset$.\\

En fait on vient de montrer la condition $Leff(X\setminus Y,X\cap H\setminus Y)$.\\

{\bf Remarque:} L'\'enonc\'e ne d\'epend pas de $Y$ tout entier mais de $Y\cap H$. On peut donc choisir $Y$ de fa\c{c}on convenable.\\

\section{\bf Cons\'equences pour le groupe de Picard}

Pour le groupe de Picard il y a des r\'esultats en plus:

\begin{theoreme} Supposons que $prof\,{\cal O}_{X\setminus Y}\ge 2$, $prof\,{\cal O}_{X\cap H\setminus Y}\ge 2$, ${\cal O}_{X\setminus Y}$ et ${\cal O}_{X\cap H\setminus Y}$ satisfont \`a la condition $F_2$ et que ${\cal O}_{X,x}$ est parafactoriel (voir [G2] XI D\'ef. 3.2) au points $x$ de
$X\setminus Y$ tels que $\dim\,\overline{\{x\}}\le \dim Y\cap H+1$. Alors $Pic\,X\setminus Y\simeq \lim\limits_{\leftarrow} Pic\,(X_n\setminus Y_n)$.
\end{theoreme}

{\bf D\'emonstration:} \`A cause de Th\'eor\`eme 4.2b) l'application $Pic\,X\setminus Y\longrightarrow \lim\limits_{\leftarrow} Pic\,X_n\setminus Y_n$ est injective. \\
Comme dans la d\'emonstration de Th\'eor\`eme 1.3, on montre qu'un \'el\'ement de $\lim\limits_{\leftarrow} Pic\,X_n\setminus Y_n$ provient d'un \'el\'ement de $Pic\,U$, $U$ \'etant un voisinage ouvert de $X\cap H\setminus Y$ dans $X\setminus Y$ convenable. Il suffit de montrer que $Pic\,X\setminus Y\longrightarrow Pic\,U$ est surjectif. Mais pour $Y':=X\setminus U$ nous savons que $\dim\, Y'\le \dim\,Y\cap H+1$, donc ${\cal O}_{X,x}$ est parafactoriel pour tout $x\in Y'\setminus Y$. Ceci implique que $(X\setminus Y,U)$ est parafactoriel, voir [G2] XI Prop. 3.3.

\begin{theoreme} Soit $X\subset \mathbb{P}_m(k)$ un sch\'ema projectif sur $k$, $Y$ un sous-espace alg\'ebrique ferm\'e, $H$ un hyperplan tel que la fl\`eche naturelle ${\cal I}\otimes{\cal O_X}\longrightarrow {\cal IO}_X$ soit injective. Soit $X_n$ le sous-espace de $X$ d\'efini par ${\cal I}^n$, $\cal I$ \'etant le faisceau d'id\'eaux qui d\'efinit $X\cap H$ dans $X$. Supposons que $prof\,{\cal O}_{X\cap H\setminus Y}\ge 3$ et que ${\cal O}_{X\cap H\setminus Y}$ satisfait $F_3$. Alors pour $n\gg 0$, on a $Pic\,(\hat{X}\setminus\hat{Y})\simeq Pic\,(X_n \setminus Y_n)$.
\end{theoreme}

{\bf D\'emonstration:} D'apr\`es [G2] XI (1.1) on a une suite exacte:
$$0\longrightarrow {\cal I}^{n+1}/{\cal I}^{n+2}\longrightarrow {\cal O}^*_{X_{n+1}}\longrightarrow {\cal O}^*_{X_n}\longrightarrow 1$$
On v\'erifie que $H^k(X\cap H\setminus Y,{\cal I}^{n+1}/{\cal I}^{n+2})=0, k=1,2,$ donc $Pic(X_{n+1}\setminus Y_{n+1})\simeq Pic(X_n\setminus Y_n), n\gg 0$, \`a cause de Th\'eor\`eme 1.2, donc $\lim\limits_\leftarrow Pic(X_n\setminus Y_n)\simeq Pic(X_n\setminus Y_n)$, $n\gg 0$.\\

{\bf D\'emonstration de Th\'eor\`eme 1.5:} On proc\`ede par r\'ecurrence sur la codimension de $X$: Le cas $X=\mathbb{P}_N(k)$ est clair parce que dans ce cas $Pic\,X=Cl\,X\simeq Cl\,X\setminus Y=Pic\,X\setminus Y$ : $X$ est lisse et $Y$ est de codimension $\ge 2$, voir [G2] XI Cor. 3.8. Nous pouvons donc supposer que $X$ n'est pas un espace projectif. Soit $X=X'\cap H$, $H$ \'etant une hypersurface et $\dim X'=\dim X+1$. Alors $X'\setminus Y$ est une intersection compl\`ete dans $\mathbb{P}_N(k)\setminus Y$.
On sait par [G2] XI Th\'eor\`eme 3.13 que ${\cal O}_{X',x}$ est parafactoriel pour tout $x\in X'\setminus Y$ avec $\dim\,\overline{\{x\}}\le \dim Y+1$ parce que $(X',x)$ est une intersection compl\`ete de dimension $\ge 4$: la codimension de $Y$ dans $X'$ est au moins $5$. Par l'hypoth\`ese de r\'ecurrence nous savons que $Pic\,X'\setminus Y\simeq \mathbb{Z}$. Par Th\'eor\`eme 5.1 (o\`u on peut remplacer ``hyperplan'' par ``hypersurface'') on obtient que $\mathbb{Z}\simeq Pic\,X'\setminus Y\simeq \lim\limits_{\leftarrow}Pic\,X'_n\setminus Y_n$. La d\'emonstration du th\'eor\`eme pr\'ec\'edant montre qu'il suffit de v\'erifier que $H^j(X\setminus Y,{\cal I}^n/{\cal I}^{n+1})=0$, $n\ge 1, j=1,2$, ${\cal I}$ \'etant le faisceau d'id\'eaux de $X$ dans $X'$. Mais il y a une r\'esolution de $(i_*({\cal I}^n/{\cal I}^{n+1}))|\mathbb{P}_N\setminus Y$, o\`u $i:X\to \mathbb{P}_N$ est l'inclusion, par des faisceaux qui sont sommes directes de ${\cal O}_{{\mathbb P}_N\setminus Y}(l)$ avec $l<0$, le nombre des faisceaux \'etant $\le k+2$ avec $k:=N-\dim\,X$. Il suffit donc de v\'erifier que $H^j(\mathbb{P}_n\setminus Y,{\cal O}_{\mathbb{P}_N}(l))=0$, $j\le k+2$, ce qui provient du fait que ${\cal H}^j_Y({\cal O}_{\mathbb{P}_N})=0$, donc $H^j_Y(\mathbb{P}_N,{\cal O}_{\mathbb{P}_N}(l))=0$, $j\le k+3$, et $H^j(\mathbb{P}_N,{\cal O}_{\mathbb{P}_N}(l))=0$, $l<0,j\le k+2$.

\section{Cas analytique}

Consid\'erons le cas $k=\mathbb{C}$. Notons que la notion de profondeur ne change pas dans le cadre analytique.\\

{\bf a) Faisceaux coh\'erents sur $X$}

\vspace{2mm}
On a un analogue de Proposition 2.1 dans le cas o\`u il y a une extension coh\'erente - ce qui n'est pas automatique dans le cas analytique, contraire au cas alg\'ebrique. 

\begin{proposition} Soit $X$ un {\sl espace analytique complexe}, $Y$ un sous-espace analytique ferm\'e, $\cal S$ un faisceau analytique coh\'erent sur $X$ (!). Supposons que $X$ est compact, $\dim Y\cap\overline{S_{l+m}({\cal S}|X\setminus Y)}<l$ pour tout $l\le \dim Y$. Alors $\dim H^s(X\setminus Y,{\cal S})<\infty$, $s<m$.
\end{proposition}

{\bf D\'emonstration:} Par [T2] III Th. 2.1, [Si2] o\`u [BS] II Th. 4.1 on sait que $R^sj_*({\cal S}|X\setminus Y)$ est coh\'erent, $s<m$, o\`u $j:X\setminus Y\longrightarrow X$. On conclut par la suite spectrale $E_2^{pq}=H^p(X,R^qj_*{\cal S})\Rightarrow H^{p+q}(X\setminus Y,{\cal S})$.

\begin{proposition} Soit $X$ une vari\'et\'e alg\'ebrique complexe compl\`ete, $Y$ une sous-vari\'et\'e ferm\'ee, $\cal S$ un faisceau alg\'ebrique coh\'erent sur $X\setminus Y$. Supposons que $\dim Y\cap\overline{S_{l+m}({\cal S})}<l$ por tout $l\le\dim Y$. Alors $(R^sj_*{\cal S})^{an}\simeq R^sj_*^{an}{\cal S}^{an}$, donc $H^s(X\setminus Y,{\cal S})\simeq H^s(X^{an}\setminus Y^{an},{\cal S}^{an})$, $s<m$.
\end{proposition}

{\bf D\'emonstration:} Pour le premier \'enonc\'e voir [Si2]. Pour le deuxi\`eme on conclut par GAGA, en comparant la suite spectrale dans la d\'emonstration pr\'ec\'edante avec l'analogue alg\'ebrique. Voir Lemme 10.8.

\begin{theoreme} Soit $X$ un sous-espace analytique ferm\'e d'un espace projectif complexe, $Y$ un sous-espace analytique ferm\'e, $\cal S$ un faisceau analytique coh\'erent sur $X$, $\cal L$ un faisceau ample sur $X$, $prof\,{\cal S}|X\setminus Y\ge m$, $\dim\,Y\cap \overline{S_{l+m}({\cal S}|X\setminus Y)}<l$ pour $l\le\dim Y$. Alors $H^s(X\setminus Y,{\cal S}\otimes_{{\cal O}_X}{\cal L}^{-\nu})=0$, $s<m$, $\nu\gg 0$.
\end{theoreme}

{\bf D\'emonstration:} Par Proposition 6.2 on peut passer au cadre alg\'ebrique et appliquer Th\'eor\`eme 1.2. Notons que $\cal S$ provient d'un faisceau alg\'ebrique coh\'erent par GAGA.\\

{\bf b) Faisceaux coh\'erents sur $X\setminus Y$}

\vspace{2mm}
Qu'est-ce qui se passe quand on ne suppose plus que $\cal S$ s'etend \`a $X$?

\vspace{2mm} Soit $X$ un {\sl espace analytique complexe}, $Y$ un sous-espace analytique ferm\'e.

\begin{lemme} {[FG]} Soit $\cal S$ un faiscau analytique coh\'erent sur $X\setminus Y$, $d\ge 0$, $\dim\,Y\le d$, $\dim S_{l+2}({\cal S})\le l$ pour tout $l\le d$. Alors il y a une extension coh\'erente $\hat{\cal S}$ de $\cal S$ \`a $X$, \`a savoir $j_*{\cal S}$, o\`u $j:X\setminus Y\longrightarrow X$ est l'inclusion, telle que $\dim\,S_{l+2}(\hat{\cal S})\le l$ pour $l\le d$.
\end{lemme}   

On a une variante de Th\'eor\`eme 6.3:

\begin{theoreme} Soit $X$ un sous-espace analytique ferm\'e d'un espace projectif complexe, $Y$ un sous-espace analytique ferm\'e, $\cal S$ un faisceau analytique coh\'erent sur $X\setminus Y$, $\cal L$ un faisceau ample sur $X$, $prof\,{\cal S}\ge m$, $\dim\,S_{l+m}({\cal S})\le l$ pour $l\le\dim Y$. Alors $H^s(X\setminus Y,{\cal S}\otimes_{{\cal O}_{X\setminus Y}}{\cal L}^{-\nu}|X\setminus Y)=0$, $s<m$, $\nu\gg 0$.
\end{theoreme}

{\bf D\'emonstration:} Il suffit de traiter le cas o\`u $prof\,{\cal S}\ge m+\dim Y+1$: On a $\dim S_{m+\dim Y}({\cal S})\le\dim Y$. Posons $Y':=Y\cup S_{m+\dim Y}({\cal S})$; alors $prof{\cal S}|(X\setminus Y')\ge m+\dim Y+1=m+\dim Y'+1$, et $H^s(X\setminus Y,{\cal S})\longrightarrow H^s(X\setminus Y',{\cal S})$ est injectif, $s<m$, \`a cause de [BS] II Th. 3.6.\\
Soit d'abord $m\ge 2$. Par le lemme pr\'ec\'edant on peut appliquer Th\'eor\`eme 6.3.\\
Ou bien on applique [HL2] Th\'eor\`eme 4.2.\\
Il reste \`a traiter le cas $m=1$ o\`u il faut d\'emontrer le lemme suivant:

\begin{lemme} Soit $X$ un sous-espace analytique ferm\'e d'un espace projectif complexe, $Y$ un sous-espace analytique ferm\'e, $\cal S$ un faisceau analytique coh\'erent sur $X\setminus Y$, $\cal L$ un faisceau ample sur $X$, ${\rm prof}\,{\cal S}\ge n:=\max\{1,\dim Y+2\}$. Soit $\cal L$ un faisceau ample sur $X$. Alors $H^0(X\setminus Y,{\cal S}\otimes_{{\cal O}_{X\setminus Y}}{\cal L}^{-\nu}|X\setminus Y)=0$, $\nu\gg 0$.
\end{lemme}

{\bf D\'emonstration:} R\'ecurrence sur $\dim Y$. On peut supposer que $X$ est un espace projectif.\par
Le cas $Y=\emptyset$ est assur\'e par [BS] IV Cor. 3.3.\par
Soit donc $Y\neq\emptyset$. On peut supposer ${\cal S}\neq 0$ et ${\cal L}={\cal O}_X(1)$. Soit $x\in Supp\,{\cal S}$. Pour presque toute forme lin\'eaire $f$ avec $f(x)=0$ la multiplication $f\cdot\ldots:{\cal S}_x\longrightarrow{\cal S}(1)_x$ est injective:\par
Soit $g$ une fonction lin\'eaire avec $g(x)\neq 0$.
En prenant des g\'en\'erateurs de ${\cal S}_x$ on trouve une suite $0={\cal S}^0_x\subset {\cal S}^1_x\subset\ldots\subset{\cal S}^r_x={\cal S}_x$ o\`u ${\cal S}^j_x/{\cal S}^{j-1}_x$ est engendr\'e par un \'el\'ement $[s_j]\neq 0$. Supposons que $f/g$ n'est pas contenu dans l'union des $Ann({\cal S}^j_x/{\cal S}^{j-1}_x)$, donc $f\cdot[s_j]\neq 0$. On v\'erifie par r\'ecurrence que la multiplication $f\cdot\ldots:{\cal S}^j_x\longrightarrow {\cal S}^j_x\otimes {\cal O}_X(1)_x$ est injective si $f$ appartient \`a un certain ouvert de Zariski non-vide. L'injectivit\'e est gard\'ee pour un voisinage de $x$ convenable.\par
Par un recouvrement d\'enombrable de $Supp\,{\cal S}$ on obtient qu'on peut trouver une fonction lin\'eaire $f$ telle que $f\cdot\ldots:{\cal S}\longrightarrow {\cal S}(1)$ est injective. Soit ${\cal I}$ l'id\'eal qui d\'efinit $f=0$; alors ${\cal L}^{-1}\simeq {\cal I}$ et ${\cal S}\otimes {\cal L}^{-l}\simeq {\cal I}^l\otimes {\cal S}$. On d\'eduit que ${\cal S}\otimes ({\cal I}^l/{\cal I}^m)\simeq {\cal I}^l\otimes{\cal S}/{\cal I}^m\otimes{\cal S}$, $m\ge l$.\par
L'hypoth\`ese de r\'ecurrence donne que
$H^0(X\setminus Y,({\cal I}^l\otimes{\cal S})/({\cal I}^{l+1}\otimes{\cal S}))=H^0(X\cap H\setminus Y,{\cal S}\otimes({\cal I}/{\cal I}^2)^l)=0$ pour $l\gg 0$ parce que $({\cal I}/{\cal I}^2)^{-1}|H$ est ample sur $X\cap H$. Donc 
$H^0(X\setminus Y,{\cal I}^{l+1}\otimes{\cal S})\longrightarrow H^0(X\setminus Y,{\cal I}^l\otimes{\cal S})$ est bijectif, et $H^0(X\setminus Y,{\cal I}^m\otimes{\cal S})\longrightarrow H^0(X\setminus Y,{\cal I}^l\otimes{\cal S})$ aussi, $l\gg 0$, $m>l$.\par
Supposons que $H^0(X\setminus Y,{\cal I}^l\otimes{\cal S})\neq 0, l\gg 0$: Soit $s$ un \'el\'ement diff\'erent de $0$. Le support de $s$ n'est pas un ensemble fini: autrement on obtient dans un point $x$ non-isol\'e du support une contribution non-triviale \`a ${\cal H}^0_x({\cal S}_x)$ en contradiction avec l'hypoth\`ese ${\rm prof}\,{\cal S}\ge 1$. Le support de $s$ coupe donc $H:f=0$.\par
Soit $x\in X\cap H\setminus Y$ avec $s_x\neq 0$. Alors $s_x\in\bigcap{\cal I}^m_x{\cal S}_x$. L'id\'eal ${\cal I}_x$ est un id\'eal propre parce que $x\in H$. ${\cal I}_x$ est donc contenu dans un id\'eal maximal. Le th\'eor\`eme d'intersection de Krull ([KK] 23.A.5) implique que $s_x=0$, ce qui est une contradiction.\\

{\bf c) Faisceaux formels}

\vspace{2mm}
Ici prenons les notations de Th\'eor\`eme 1.6. 

\vspace{2mm}
Pour les faisceaux formels, on peut montrer

\begin{proposition}
Supposons que ${\cal S}$ un faisceau alg\'ebrique coh\'erent sur $\hat{X}$, ${\cal I}\otimes{\cal S}|\hat{X}\setminus \hat{Y}\simeq{\cal IS}|\hat{X}\setminus \hat{Y}$, $prof\,{\cal S}/{\cal I}{\cal S}\ge m$, ${\cal S}/{\cal I}{\cal S}$ satisfait $(F_m)$. Alors

$$\begin{array}{ccc}
H^k(\hat{X}\setminus\hat{Y},{\cal S})&\stackrel{\simeq}{\longrightarrow}&H^k(X_n\setminus Y_n,{\cal S}/{\cal I}^n{\cal S})\\
\downarrow\simeq&&\downarrow\simeq\\
H^k(\hat{X}^{an}\setminus\hat{Y}^{an},{\cal S}^{an})&\stackrel{\simeq}{\longrightarrow}&H^k(X_n^{an}\setminus Y_n^{an},{\cal S}^{an}/({\cal I}^{an})^n{\cal S}^{an})
\end{array}$$
pour $k<m, n\gg 0$.
\end{proposition}

{\bf D\'emonstration:} L'isomorphisme sup\'erieur est assur\'e par Proposition 3.5, celui \`a droite par Proposition 6.2. Par Th\'eor\`eme 6.3 on obtient pour $n\gg 0, k<m$: $H^k(\hat{X}^{an}\setminus \hat{Y}^{an}, ({\cal I}^n{\cal S}/{\cal I}^{n+1}{\cal S})^{an})=0, n\gg 0$, donc $H^k(\hat{X}^{an}\setminus \hat{Y}^{an},{\cal S}^{an}/({\cal I}^{n+1}{\cal S})^{an})\longrightarrow H^k(\hat{X}^{an}\setminus \hat{Y}^{an}, {\cal S}^{an}/({\cal I}^n{\cal S})^{an})$ est injectif, donc bijectif. Par cons\'equent, $H^k(X_n^{an}\setminus Y_n^{an},{\cal S}^{an}/({\cal I}^{an})^n{\cal S}^{an})\simeq \lim\limits_\leftarrow H^k(X_m^{an}\setminus Y_m^{an},{\cal S}^{an}/({\cal I}^{an})^m{\cal S}^{an})\simeq H^k(\hat{X}^{an}\setminus \hat{Y}^{an},{\cal S}^{an})$, voir d\'emonstration de Proposition 3.5. La fl\`eche inf\'erieure est donc bijective.

\begin{proposition}
Soit $\cal S$ un faisceau alg\'ebrique coh\'erent sur $X\setminus Y$ tel que ${\cal I}\otimes{\cal S}\simeq{\cal IS}$, $prof\,{\cal S}\ge m$, $\dim\,Y\cap\overline{S_{l+m}({\cal S})}< l$ pour $l\le\dim\,Y$. Alors, pour $k<m-1$:
$$\begin{array}{ccc}
H^k(X\setminus Y,{\cal S})&\simeq&H^k(\hat{X}\setminus\hat{Y},\hat{\cal S})\\
\downarrow\simeq&&\downarrow\simeq\\
H^k(X^{an}\setminus Y^{an},{\cal S}^{an})&\simeq&H^k(\hat{X}^{an}\setminus\hat{Y}^{an},\hat{\cal S}^{an})
\end{array}$$
\end{proposition}

{\bf D\'emonstration:} L'isomorphie \`a gauche d\'ecoule de Proposition 6.2, celle d'en haut de Th\'eor\`eme 3.1, celle \`a droite de Proposition 6.7. L'isomorphie en bas en r\'esulte mais s'obtient aussi comme dans Th\'eor\`eme 3.1, en utilisant Th\'eor\`eme 6.3 au lieu de 1.2; notons que $\cal S$ admet une extension coh\'erente \`a $X$ .

\begin{proposition}
Soit ${\cal S}$ un faisceau analytique coh\'erent sur $\hat{X}^{an}\setminus\hat{Y}^{an}$, ${\cal I}^{an}\otimes{\cal S}\simeq{\cal I}^{an}{\cal S}$, $\dim S_{l+2}({\cal S}/{\cal I}^{an}{\cal S})\le l$, $l\le\dim Y\cap H$. Alors $\cal S$ provient d'un faisceau alg\'ebrique coh\'erent sur $\hat{X}\setminus\hat{Y}$. 
\end{proposition}

{\bf D\'emonstration:} \`A cause de Lemme 6.4 il y a pour tout $n$ une extension coh\'erente de ${\cal S}/({\cal I}^{an})^n{\cal S}$ \`a $X_n^{an}$. Par GAGA la derni\`ere est alg\'ebrique, il y a donc un faisceau alg\'ebrique ${\cal T}_n$ sur $X_n\setminus Y_n$ tel que ${\cal T}_n^{an}\simeq {\cal S}/({\cal I}^{an})^n{\cal S}$. On a donc $({\cal T}_{n+1}/{\cal I}^n{\cal T}_{n+1})^{an}\simeq {\cal T}^{an}_{n+1}/({\cal I}^{an})^n{\cal T}_{n+1}^{an}\simeq {\cal S}/({\cal I}^{an})^n{\cal S}\simeq {\cal T}_n^{an}$. Par cons\'equent, $j_*^{an}({\cal T}_{n+1}/{\cal I}^n{\cal T}_{n+1})^{an}\simeq j_*^{an}{\cal T}_n^{an}$, c.-\`a d. $(j_*({\cal T}_{n+1}/{\cal I}^n{\cal T}_{n+1}))^{an}\simeq (j_*{\cal T}_n)^{an}$, voir Proposition 6.2. Par GAGA, $j_*({\cal T}_{n+1}/{\cal I}^n{\cal T}_{n+1})\simeq j_*{\cal T}_n$, d'o\`u ${\cal T}_{n+1}/{\cal I}^n{\cal T}_{n+1}\simeq {\cal T}_n$. Soit ${\cal T}:=\lim\limits_\leftarrow {\cal T}_n$, alors ${\cal T}$ est coh\'erent par [H1] II Prop. 9.6, et ${\cal T}^{an}\simeq {\cal S}$.

\begin{proposition} Soit ${\cal S}$ un faisceau analytique coh\'erent sur $\hat{X}^{an}\setminus\hat{Y}^{an}$, ${\cal I}^{an}\otimes{\cal S}\simeq{\cal I}^{an}{\cal S}$, $m\ge 2$, $\dim S_{l+m}({\cal S}/{\cal I}^{an}{\cal S})\le l$ pour $l\le\dim Y\cap H$. Alors
$H^k(\hat{X}^{an}\setminus\hat{Y}^{an},{\cal S})\simeq H^k(X_n^{an}\setminus Y_n^{an},{\cal S}/({\cal I}^{an})^n{\cal S})$ pour $k<m, n\gg 0$.
\end{proposition}

{\bf D\'emonstration:} On proc\`ede comme dans la d\'emonstation de Proposition 6.8. Au lieu de Th\'eor\`eme 6.3 on applique Th\'eor\`eme 6.5. La finitude de $H^k(X_n^{an}\setminus Y_n^{an},{\cal S}/({\cal I}^{an})^n{\cal S})$ provient de Lemme 6.4 et Proposition 6.1.

\section{Fibr\'es vectoriels analytiques}

Soit maintenant $X$ un {\sl sous-espace analytique} ferm\'e de $\mathbb{P}_N(\mathbb{C})$ et $Y$ un sous-espace analytique ferm\'e de $X$. Par Chow $X$ et $Y$ sont alg\'ebriques, c.-\`a d. proviennent de sch\'emas complexes. Soit $H$ un hyperplan, $d:=\dim\,Y\cap H$, o\`u $\dim\emptyset:=-1$.\\

{\bf a) Sections de faisceaux}

\begin{theoreme}
Soit $\cal S$ un faisceau analytique coh\'erent sur $X$ tel que $\dim\,S_{l+2}({\cal S})\le l$ pour $l\le d+1$. Soit $s\in\Gamma(X\cap H\setminus Y,{\cal S})$. Alors $s$ admet une extension unique \`a une section $\hat{s}\in\Gamma(X,{\cal S})$.
\end{theoreme}

{\bf D\'emonstration:} R\'ecurrence sur $d$. Cas $d=-1$: D'apr\`es [Go] II Th. 3.3.1, on peut \'etendre $s$ \`a une section dans $\Gamma(U,{\cal S})$, $U$ \'etant un voisinage de $X\cap H$ dans $X$ dont le compl\'ementaire est un ferm\'e de Stein. Or, $prof\,{\cal S}\ge 2$, donc $\Gamma(X,{\cal S})\simeq \Gamma(U,{\cal S})$ parce que $H^j_{X\setminus U}(X,{\cal S})=0, j\le 1$, par [BS] I Theorem 3.1.\\
Passage de $d-1$ \`a $d$: Nous pouvons supposer que $X=\mathbb{P}_N$ et $Y=Y\cap H$ est r\'eduit.\\
Soit $Y':=Sing\;Y$, donc $\dim\,Y'\le d-1$.\\
Fixons un point de $Y\setminus Y'$. Choisissons des coordonn\'ees locales $z_1,\ldots,z_N$ tel que notre point correspond \`a $0$ et qu'on a un voisinage de la forme $U:=\{|z_j|<\epsilon,j=1,\ldots, N\}$. Nous pouvons achever que $H\cap U=\{z\in U\,|\,z_N=0\}$ et $Y\cap U=\{z\in U\,|\,z_{d+1}=\ldots=z_N=0\}$. \\
Soient $(\epsilon_\nu)_{\nu\ge 1}$ et $(\delta_\nu)_{\nu\ge 1}$ deux suites strictement d\'ecroissantes avec $\epsilon>\epsilon_\nu\searrow 0$, $\epsilon>\delta_\nu\searrow 0$. Soit $W_\nu:=\{z\in U\,|\,\max(|z_{d+1}|,\ldots,|z_{N-1}|)>\epsilon_\nu, |z_N|<\delta_\nu\}$, $W:=\bigcup W_\nu$. Les fermetures des $W$ dans $\overline{U}\setminus Y$ forment un syst\`eme fondemental de voisinages ferm\'es de $\overline{U}\cap H\setminus Y$ dans $\overline{U}\setminus Y$.\\
Or, $s|U\cap H\setminus Y$ s'\'etend \`a une section $s'$ sur un tel $W$. Nous voulons montrer qu'il y a une extension unique $s''$ de $s'$ sur $U_1\setminus Y$, avec $U_1:=U\cap\{|z_N|<\delta_1\}$.\\
Soit $V_\nu:=W\cup\{z\in U\,|\,\delta_\nu<|z_N|<\delta_1$ ou $\max(|z_{d+1}|,\ldots,|z_{N-1}|)>\epsilon_\nu,|z_N|<\delta_1\}$. Montrons par r\'ecurrence que $s'$ admet une extension unique \`a $s_\nu\in\Gamma(V_\nu,{\cal S})$.\\
Passage de $\nu$ \`a $\nu+1$:\\
D'apr\`es [Si3] Proposition 3.14, p. 141 (en y posant $\cup$ au lieu de $\cap$; voir aussi [Si3] Prop. 3.13, p. 140) il y a une extension unique $t$ de $s_\nu|U\cap\{\max(|z_{d+1}|,\ldots,|z_{N-1}|)>\epsilon_\nu,\delta_{\nu+2}<|z_N|<\delta_1\}$ \`a $U\cap\{\delta_{\nu+2}<|z_N|<\delta_1\}$. La restriction de $t$ \`a $U\cap\{\delta_{\nu}<|z_N|<\delta_1\}$ doit coincider avec celle de $s_\nu$. En plus, la restriction de $t$ \`a $U\cap\{\delta_{\nu+2}<|z_N|<\delta_{\nu+1}\}$ doit coincider avec l'extension de $s_\nu|U\cap\{\delta_{\nu+2}<|z_N|<\delta_{\nu+1},\max(|z_{d+1}|,\ldots,|z_{N-1}|) >\epsilon_{\nu+1}\}$. En total $t$ est une extension de $s_\nu|V_\nu\cap\{|z_N|>\delta_{\nu+2}\}$. En utilisant $t$ on obtient donc l'extension $s_{\nu+1}$.\\
Les $s_\nu$ conduisent \`a une extension unique de $s'$ \`a $U_1\setminus Y$. D'apr\`es [BS] II Theorem 3.6 on obtient une extension unique \`a $U_1$.\\
Le germe de l'extension en $0$ ne d\'epend que de $s$: Supposons que l'on a deux extensions. Elles doivent coincider sur un ensemble $W$ de la forme indiqu\'ee en haut. Gr\^ace \`a l'unicit\'e de l'extension d'une section sur $W$ on obtient l'\'enonc\'e d'unicit\'e.\\
On obtient donc une extension de $s$ \`a $X\setminus Y'$, donc par r\'ecurrence \`a $X$.\\
Le proc\'ed\'e ici est semblable \'a celui dans le cas local dans [Ha2].

\begin{corollaire}
Soit $\cal S$ un faisceau analytique coh\'erent sur $X\setminus Y$ tel que $\dim\,S_{l+2}({\cal S})\le l$ pour $l\le d+1$. Soit $s\in\Gamma(X\cap H\setminus Y,{\cal S})$. Alors $s$ admet une extension unique \`a une section $\hat{s}\in\Gamma(X\setminus Y,{\cal S})$.
\end{corollaire}

{\bf D\'emonstration:} D'apr\`es Lemme 6.4, $j_*{\cal S}$ est coh\'erent et satisfait \`a l'hypoth\`ese de Th\'eor\`eme 7.1. On a une extension unique \`a une section de $j_*{\cal S}$, il y a donc une extension \`a une section de $\cal S$ qui est unique aussi \`a cause de [BS] II Theorem 3.6.\\

{\bf b) Extension de faisceaux}

\begin{theoreme} Soit $\cal S$ un $({\cal O}_{X\setminus Y}|X\cap H\setminus Y)$-module coh\'erent, $d:=\dim\,Y\cap H$, $\dim S_{l+2}({\cal S})\le l$ pour $l\le d+1$.\hfill(*)\\ 
Alors il y a un $({\cal O}_X|X\cap H)$-module coh\'erent $\hat{\cal S}$ tel que $\hat{\cal S}|X\cap H\setminus Y={\cal S}$ et $\dim\;S_{l+2}(\hat{S})\le l$ pour $l\le d+1$, \`a savoir $j_*{\cal S}$, o\`u $j:X\cap H\setminus Y\to X\cap H$ est l'inclusion. 
\end{theoreme}

{\bf D\'emonstration:} R\'ecurrence sur $d$. Le cas $d=-1$, c.-\`a d. $Y\cap H=\emptyset$ est trivial.\\
Passage de $d-1$ \`a $d$: Nous pouvons supposer que $X=\mathbb{P}_N$ et $Y=Y\cap H$ est r\'eduit.\\
Soit $Y':=Sing\;Y$, donc $\dim\,Y'\le d-1$.\\
Fixons un point de $Y\setminus Y'$. Choisissons coordonn\'ees locales $z_1,\ldots,z_N$ tel que notre point correspond \`a $0$ et qu'on a un voisinage de la forme $U:=\{|z_j|<\epsilon,j=1,\ldots, N\}$. Nous pouvons achever que $H\cap U=\{z\in U\,|\,z_N=0\}$ et $Y\cap U=\{z\in U\,|\,z_{d+1}=\ldots=z_N=0\}$. \\
Soient $(\epsilon_\nu),(\delta_\nu)$ sont des suites strictement d\'ecroissantes avec $\epsilon>\epsilon_\nu\searrow 0, \epsilon>\delta_\nu\searrow 0$.
Posons $W_\nu:=\{z\in U\,|\,\max(|z_{d+1}|,\ldots,|z_{N-1}|)>\epsilon,|z_N|<\delta\}$ et $W=\bigcup W_\nu$. Les fermetures des tels $W$ dans $\overline{U}\setminus Y$ forment un syst\`eme fondemental de voisinages ferm\'es de $\overline{U}\cap H\setminus Y$ dans $\overline{U}\setminus Y$.\\
D'apr\`es Corollaire 10.2 le faisceau ${\cal S}|U\cap H\setminus Y$ est repr\'esent\'e par un faisceau coh\'erent ${\cal S}_1$ sur un tel $W$. Nous pouvons supposer que les conditions sur la profondeur (*) du th\'eor\`eme sont gard\'ees.\\
Nous voulons montrer qu'il y a une extension ${\cal S}_2$ de ${\cal S}_1$ sur $U_1\setminus Y$ avec les conditions (*) de profondeur, ou $U_1:=U\cap\{|z_N|<\delta_1\}$\\
On proc\`ede ici comme dans la d\'emonstration de Th\'eor\`eme 7.1, en utilisant [Si3] Theorem 7.4, p. 243. Remarquons qu'on a unicit\'e \`a isomorphisme pr\`es mais ceci suffit pour notre but.\\
Maintenant il y a une extension coh\'erente ${\cal S}_3$ de ${\cal S}_2$ \`a $U_1$, o\`u l'on garde les conditions (*), par Lemme 6.4.\\
Or, $H^0(U_1\cap\{|z_N|<\delta_1\},{\cal S}_3)\simeq H^0(W,{\cal S}_1)$: On montre d'abord que $H^0(U_1\cap\{|z_N|<\delta_1\}\setminus Y,{\cal S}_2)\simeq H^0(W,{\cal S}_1)$. Ici on utilise [Si3] Proposition 3.14, p. 141, comme dans la d\'emonstration de Th\'eor\`eme 7.1. Il reste \`a montrer que $H^0(U_1,{\cal S}_3)\simeq H^0(U_1\setminus Y,{\cal S}_2)$, o\`u on utilise [BS] II Theorem 3.6.\\
En particulier ceci montre que ${\cal S}_3\simeq (j_W)_*{\cal S}_1$, o\`u $j_W:W\setminus U_1$ est l'inclusion, donc ${\cal S}_3|U\cap H$ est l'image directe de ${\cal S}|U\cap H\setminus Y$ par rapport \`a l'inclusion.\\
On obtient donc que $j_*{\cal S}|X\cap H\setminus Y'$ est un $({\cal O}_X|X\cap H\setminus Y')$-module coh\'erent avec (*). Par r\'ecurrence on conclut que $j_*{\cal S}$ a les propri\'et\'es d\'esir\'ees.

\begin{theoreme} Sous les hypoth\`eses du th\'eor\`eme pr\'ec\'edant il y a une extension de $\cal S$ \`a un faisceau analytique coh\'erent $\hat{\cal S}$ sur $X$ tel que $\dim\;S_{l+2}(\hat{S})\le l$ pour $l\le d+1$. L'extension est unique \`a isomorphie pr\`es.
\end{theoreme}

{\bf D\'emonstration:} Existence: D'apr\`es le th\'eor\`eme pr\'ec\'edant nous pouvons supposer que $Y\cap H=\emptyset$, donc $Y$ consiste d'un nombre fini de points. De nouveau nous pouvons nous borner au cas $X=\mathbb{P}_N$. D'apr\`es Cor. 10.2 il y a une extension de $\cal S$ \`a un faisceau analytique coh\'erent ${\cal S}_1$ sur un voisinage $W$ de $H$ dans $X$. On peut supposer que $W=H\cup\{z\,|\,\max(|z_1|,\ldots,|z_N|)>R\}$, o\`u on a identifi\'e $X\setminus H$ avec $\mathbb{C}^N$. D'apr\`es [Si3] Theorem 7.4 il y a une extension analytique coh\'erente de ${\cal S}_1|W\setminus H$ \`a $X\setminus H$. donc de ${\cal S}_1$ \`a $X$.\\
Unicit\'e: Soient ${\cal S}_1$ et ${\cal S}_2$ deux telles extensions. L'identit\'e donne une section de ${\cal T}:=Hom({\cal S}_1,{\cal S}_2)$ sur $X\cap H\setminus Y$. D'apr\`es [Si3] Prop. 3.13, p. 140 $\cal T$ dispose des m\^emes propri\'et\'es de profondeur que ${\cal S}_2$. Par Th\'eor\`eme 7.1 on peut donc trouver une extension de cette section \`a $X$ qui est unique. On continue comme d'habitude.

\begin{corollaire}
Sous les hypoth\`eses de Th\'eor\`eme 7.3 il y a une extension de $\cal S$ \`a un faisceau analytique coh\'erent $\hat{\cal S}$ sur $X\setminus Y$ tel que $\dim\;S_{l+2}(\hat{S})\le l$ pour $l\le d+1$. L'extension est unique \`a isomorphie pr\`es.
\end{corollaire}

{\bf D\'emonstration:} L'existence d\'ecoule directement du th\'eor\`eme pr\'ec\'edant.\\
Unicit\'e: Soient $\hat{\cal S}_1, \hat{\cal S}_2$ deux telles extensions. Par Lemme 6.4, $j_*\hat{\cal S}_1$ et $j_*\hat{\cal S}_2$ sont coh\'erents  et satisfont \`a (*), ils sont donc isomorphes, donc $\hat{\cal S}_1$ et $\hat{\cal S}_2$ aussi.\\

Les \'enonc\'es suivants pr\'eparent le prochain paragraphe:

\begin{theoreme} Soit $Y$ un sous-espace analytique ferm\'e d'un espace analytique complexe $X$, $\dim\, Y\cap S_{l+2}({\cal O}_X)\le l$ pour tout $l$. Alors $Vect\,X\longrightarrow Vect\,(X\setminus Y)$ est injectif.
\end{theoreme}

{\bf D\'emonstration:} On proc\`ede comme dans la d\'emonstration de Th\'eor\`eme 4.1, en utilisant [BS] II Th. 3.6.

\begin{proposition}
Soit $prof\,{\cal O}_{X\setminus Y}\ge 2$, $\dim S_{l+2}({\cal O}_{X\setminus Y})\le l$ pour tout $l\le\dim Y$. Alors l'application $Vect(X\setminus Y)\longrightarrow \lim\limits_\leftarrow Vect(X_n\setminus Y_n)$ est injective.
\end{proposition}

{\bf D\'emonstration:} On proc\`ede comme dans la d\'emonstration de Th\'eor\`eme 4.2b), en utilisant Th\'eor\`eme 6.5.

\section{Fibr\'es vectoriels alg\'ebriques et analytiques}

Retournons au contexte alg\'ebrique.

\begin{proposition} Soit $X$ une vari\'et\'e alg\'ebrique complexe, $X$ complet, $Y$ une sous-vari\'et\'e ferm\'ee, $\dim S_{l+2}({\cal O}_{X\setminus Y})\le l$ pour $l\le \dim Y$, alors $Vect(X\setminus Y)\simeq Vect(X^{an}\setminus Y^{an})$.
\end{proposition}

{\bf D\'emonstration:} Soit ${\cal S}$ un repr\'esentant d'un \'el\'ement de $Vect^{an}(X^{an}\setminus Y^{an})$. Alors  $\dim S_{l+2}({\cal S})\le l$ pour $l\le \dim Y$. Par Lemme 6.4, il y a une extension \`a un faisceau analytique coh\'erent sur $X^{an}$ qui est alg\'ebrique par GAGA, $\cal S$ provient donc d'un fibr\'e vectoriel alg\'ebrique.\\
Injectivit\'e: Soit $j:X\setminus Y\longrightarrow X$ l'inclusion, et soient ${\cal S}_1, {\cal S}_2$ localement libres sur $X\setminus Y$ tels que ${\cal S}_1^{an}\simeq {\cal S}_2^{an}$. Alors $(j_*{\cal S}_1)^{an}\simeq j_*^{an}{\cal S}_1^{an}\simeq j_*^{an}{\cal S}_2^{an}\simeq (j_*{\cal S}_2)^{an}$, voir [Si2], donc $j_*{\cal S}_1\simeq j_*{\cal S}_2$. Ceci implique ${\cal S}_1\simeq {\cal S}_2$.

\vspace{2mm}
Reprenons les anciennes hypoth\`eses.

\begin{proposition} Soit $\dim S_{l+2}({\cal O}_{X\cap H\setminus Y})\le l$ pour $l\le\dim Y\cap H$. Alors
$$\begin{array}{ccc}
Vect(\hat{X}\setminus\hat{Y})&\stackrel{\simeq}{\longrightarrow}&\lim\limits_\leftarrow Vect(X_n\setminus Y_n)\\
\downarrow\simeq&&\downarrow\simeq\\
Vect(\hat{X}^{an}\setminus\hat{Y}^{an})&\stackrel{\simeq}{\longrightarrow}&\lim\limits_\leftarrow Vect(X_n^{an}\setminus Y_n^{an})
\end{array}$$
\end{proposition}

{\bf D\'emonstration:} Par Th\'eor\`eme 4.7 on a l'isomorphisme sup\'erieur, par la proposition pr\'ec\'edante celui \`a droite. Il suffit donc de montrer que la fl\`eche \`a gauche est surjective ou bien que la fl\`eche inf\'erieure est injective.\\
Pour l'injectivit\'e de la fl\`eche inf\'erieure on peut proc\'eder comme dans la d\'emonstration de Th\'eor\`eme 4.7, en utilisant Proposition 6.10.

\begin{theoreme} Soit $\dim S_{l+2}({\cal O}_X)\le l$ pour $l\le \dim Y\cap H+1$ le long de $X\cap H\setminus Y$, alors
$$\lim_\rightarrow Vect(U)\simeq \lim_\rightarrow Vect(U^{an})\simeq H^1(X^{an}\cap H^{an}\setminus Y^{an}, Gl({\cal O}_{X^{an}}))$$
o\`u $U$ parcourt les voisinages ouverts (de Zariski) de $X\cap H\setminus Y$ dans $X\setminus Y$.
\end{theoreme}

{\bf D\'emonstration:} D'abord \'etudions la fl\`eche \`a droite.\\
Si $l\le d+1$, la fermeture de chaque composante irr\'eductible de $S_{l+2}({\cal O}_{X\setminus H})$ de dimension $>l$ ne rencontre pas $X\cap H\setminus Y$, sa dimension est donc $\le d+1$. En agrandissant $Y$ on peut donc supposer que $\dim S_{l+2}({\cal O}_{X\setminus Y})\le l$ pour $l\le d+1$. \\
Surjectivit\'e: Un \'el\'ement de $H^1(X^{an}\cap H^{an}\setminus Y^{an},Gl({\cal O}_{X^{an}}))$ est repr\'esent\'e par un fibr\'e vectoriel $\cal S$ sur un voisinage de $X^{an}\cap H^{an}\setminus Y^{an}$ dans $X^{an}\setminus Y^{an}$, voir Corollaire 10.4. Par hypoth\`ese, $\dim S_{l+2}({\cal S})\le l$ pour $l\le d+1$. Par Th\'eor\`eme 7.4 on peut donc supposer que $\cal S$ admet une extension coh\'erente ${\cal S}_1$ sur $X^{an}$. Soit $\Sigma$ le sous-sch\'ema de $X$ qui correspond au lieu o\`u ${\cal S}_1$ n'est pas localement libre, et soit $U:=X\setminus Y\cup\Sigma$. Alors $\Sigma\cap H=\emptyset$, et ${\cal S}_1|U^{an}$ est le fibr\'e vectoriel cherch\'e.\\
Injectivit\'e: Supposons que ${\cal S}_1, {\cal S}_2$ sont des fibr\'es vectoriels sur $U^{an}$ tels que ${\cal S}_1|X^{an}\cap H^{an}\setminus Y^{an}\simeq{\cal S}_2|X^{an}\cap H^{an}\setminus Y^{an}$. Par Corollaire 7.5 on a ${\cal S}_1\simeq{\cal S}_2$.\\
Pour la bijectivit\'e de la fl\`eche \`a gauche, il suffit de montrer que $Vect(X\setminus \Sigma)\simeq Vect(X^{an}\setminus \Sigma^{an})$ si $Y\subset \Sigma$, $\Sigma$ ferm\'e dans $X$, $\Sigma\cap H=Y\cap H$. Ceci d\'ecoule de Proposition 8.1.\\

{\bf Remarque:} L'hypoth\`ese sur ${\cal O}_X$ est assur\'ee si $X$ est normal le long de $X\cap H\setminus Y$ et $codim_XY\cap H\ge 4$.\\

{\bf D\'emonstration de Th\'eor\`eme 1.6:} D'abord l'hypoth\`ese sur $H$ implique que $prof\,{\cal O}_{X,x}=prof\,{\cal O}_{X\cap H,x}+1$ pour $x\in X\cap H\setminus Y$. On a donc pr\`es de $X\cap H\setminus Y$:\\
$\dim S_{l+2}({\cal O}_{X\setminus Y})=\dim S_{l+2}({\cal O}_{X\setminus Y})\cap H + 1=\dim S_{l+1}({\cal O}_{X\cap H\setminus Y})+1\le l$ pour $l\le\dim Y\cap H+1$. On dispose alors de l'hypoth\`ese du th\'eor\`eme pr\'ec\'edant.\\
Par Th\'eor\`eme 1.3 on a la bijectivit\'e des fl\`eches sup\'erieures. Par le th\'eor\`eme pr\'ec\'edant, la fl\`eche verticale \`a gauche et la premi\`ere fl\`eche inf\'erieure sont bijectives. En plus, $Vect(X_n\setminus Y_n)\simeq Vect(X_n^{an}\setminus Y_n^{an})$ par Proposition 8.1. Par Proposition 8.2 on conclut donc que les fl\`eches restantes sont bijectives aussi.\\

{\bf Remarque:} L'hypoth\`ese sur ${\cal O}_{X\cap H\setminus Y}$ est donn\'ee en particulier si $X\cap H\setminus Y$ est normal et $codim_XY\cap H\ge 4$.

\begin{proposition} Supposons que $\dim S_{l+2}({\cal O}_{X\cap H\setminus Y})\le l$ pour tout $l\le\dim Y\cap H$. Pour $n\gg 0$, on a un diagramme commutatif
$$\begin{array}{ccc}
\lim\limits_{\stackrel{\leftarrow}{m}}Vect\,(X_m\setminus Y_m)&\longrightarrow& Vect\,(X_n\setminus Y_n)\\
\downarrow\simeq&&\downarrow\simeq\\
\lim\limits_{\stackrel{\leftarrow}{m}}Vect(X_m^{an}\setminus Y_m^{an})&\longrightarrow& Vect(X_n^{an}\setminus Y_n^{an})
\end{array}$$
o\`u le noyau de chaque application horizontale est trivial.
\end{proposition}

{\bf D\'emonstration:} On utilise Proposition 4.5 et 8.1.\\

On peut comparer avec $Vect\,(X\cap H\setminus Y)$, en utilisant une g\'en\'eralisation du th\'eor\`eme de Kodaira, o\`u $Sing(X)$ d\'enote le lieu singulier de $X$:

\begin{theoreme} (voir [HL2] Cor. 4.4) Soit $m\in\mathbb{N}$, ${\rm codim}_XY\ge 
m+1$, $prof_{Sing(X)\setminus Y}{\cal O}_{X\setminus Y}\ge m$, $\cal L$ un faisceau ample sur $X$. 
Alors
$H^q(X\setminus Y,{\cal L}^{-1})=H^q(X^{an}\setminus Y^{an},({\cal L}^{an})^{-1})=0$ pour $q<m$.
\end{theoreme}

Alors on en d\'eduit:
\begin{theoreme} Soit ${\rm codim}_{X\cap H}Y\cap H\ge 3$, $prof_{Sing(X\cap H)\setminus Y}{\cal O}_{X\cap H\setminus Y}\ge 2$. Alors on a un diagramme commutatif
$$\begin{array}{ccc}
\lim\limits_{\stackrel{\leftarrow}{m}}Vect\,(X_m\setminus Y_m)&\longrightarrow& Vect\,(X\cap H\setminus Y)\\
\downarrow\simeq&&\downarrow\simeq\\
\lim\limits_{\stackrel{\leftarrow}{m}}Vect(X_m^{an}\setminus Y_m^{an})&\longrightarrow& Vect(X^{an}\cap H^{an}\setminus Y^{an})
\end{array}$$
o\`u le noyau des applications horizontales est trivial.
\end{theoreme}

{\bf D\'emonstration:} 
La bijectivit\'e des fl\`eches verticales d\'ecoule de Proposition 8.1.\\
Soit $X_n$ le sous-espace de $X$ d\'efini par ${\cal I}^n$. 
De fa\c{c}on analogue \`a [G2] XI (1.1), on a une suite exacte, voir Lemme 4.4:
$$0\longrightarrow (({\cal I}^{an})^k/({\cal I}^{an})^{k+1})^{\oplus n^2}\longrightarrow GL_n({\cal O}_{X^{an}_{k+1}})\longrightarrow GL_n({\cal O}_{X^{an}_k})\longrightarrow 1$$
On v\'erifie que $H^1(X\cap H\setminus Y,({\cal I}^{an})^k/({\cal I}^{an})^{k+1})=0$ pour tout $k\ge 1$, d'apr\`es Th\'eor\`eme 8.5; notons que $({\cal I}^{an})^k/({\cal I}^{an})^{k+1}\simeq ({\cal I}^{an}/({\cal I}^{an})^2)^k$.\\

{\bf D\'emonstration de Th\'eor\`eme 1.7:} D'apr\`es Th\'eor\`eme 4.2, l'application $Vect(X\setminus Y)\to\lim\limits_\leftarrow Vect(X_n\setminus Y_n)$ est injective. Th\'eor\`eme 8.5 implique que le noyau de $\lim\limits_\leftarrow Vect(X_n\setminus Y_n)\to Vect(X\cap H\setminus Y)$ est trivial. D'apr\`es Proposition 8.1 resp. Th\'eor\`eme 8.6 il est \'egal si l'on travaille dans le contexte alg\'ebrique ou analytique.\\

\section{Fibr\'es \`a connexion}

{\bf a) Cas g\'en\'eral}\\

Pour des fibr\'es \`a connexion il y a des r\'esultats plus complets.

\vspace{2mm}
Si $X$ est une vari\'et\'e alg\'ebrique complexe lisse ou bien une vari\'et\'e analytique complexe, soit $Vect_c(X)$ resp. $Vect_{ci}(X)$ l'ensemble des classes d'isomorphie de fibr\'es vectoriels sur $X$ \`a connexion resp. connexion int\'egrable, voir Appendice et [D] I.2. Dans le contexte alg\'ebrique nous consid\'erons aussi le cas des fibr\'es \`a connexion int\'egrable r\'eguli\`ere o\`u nous \'ecrivons $Vect_{cir}(X)$, voir [D] II.4.

\vspace{2mm}
D'abord

\begin{theoreme} Soit $X$ une vari\'et\'e analytique, $Y$ un sous-ensemble analytique ferm\'e de codimension $\ge 2$. Alors $Vect_c\,X\simeq Vect_c(X\setminus Y)$.
\end{theoreme}

{\bf D\'emonstration:} D'apr\`es Th\'eor\`eme 7.6 la fl\`eche $Vect\;X\to Vect(X\setminus Y)$ est injective, donc $Vect_cX\to Vect_c(X\setminus Y)$ aussi: Soit $\cal E$ coh\'erent et localement libre sur $X$. Une connexion sur ${\cal E}|X\setminus Y$ peut \^etre donn\'ee par une application ${\cal O}_{X\setminus Y}$-lin\'eaire $D_1:P^1({\cal E}|X\setminus Y)\to \Omega^1_{X\setminus Y}\otimes{\cal E}|X\setminus Y$ avec $D_1\circ i=id$, voir Proposition 10.6. Or, $D_1$ correspond \`a une section ${\rm D}_1$ de $Hom(P^1({\cal E}|X\setminus Y),\Omega^1_{X\setminus Y}\otimes{\cal E}|X\setminus Y)\simeq {\cal T}|X\setminus Y$, o\`u ${\cal T}:=Hom(P^1({\cal E}),\Omega^1_X\otimes{\cal E})$.  Or, ${\rm D}_1$ s'\'etend de fa\c{c}on unique \`a une section de $\cal T$ qui correspond \`a une connexion sur $\cal E$, d'apr\`es [BS] II Theorem 3.6. \\
Surjectivit\'e: Soit $\cal E$ coh\'erent et localement libre sur $X\setminus Y$ avec une connexion et $j:X\setminus Y\to X$ l'inclusion. Alors $j_*{\cal E}$ est coh\'erent d'apr\`es Lemme 6.4. Par un raisonnement comme dans le cas de l'injectivit\'e, la connexion s'\'etend \`a une connexion unique sur $j_*{\cal E}$ d'apr\`es [BS] II Theorem 3.6. Par cons\'equent, $j_*{\cal E}$ est localement libre, voir [Ma] Remark (1.2).\\

Dans le cas alg\'ebrique on a

\begin{theoreme} Soit $X$ une vari\'et\'e alg\'ebrique complexe, $Y$ un sous-ensemble alg\'ebrique ferm\'e de codimension $\ge 2$.\\
a) Si $X$ est lisse on a $Vect_c\,X\simeq Vect_c(X\setminus Y)$.\\
b) Si $X$ est complet et $X\setminus Y$ lisse on a $Vect_c(X\setminus Y)\simeq Vect_c(X^{an}\setminus Y^{an})$.
\end{theoreme}

{\bf D\'emonstration:} a) On raisonne comme dans la d\'emonstration pr\'ec\'edante. D'apr\`es Th\'eor\`eme 4.1 on a que la fl\`eche $Vect\;X\to Vect(X\setminus Y)$ est injective, donc $Vect_cX\to Vect_c(X\setminus Y)$ aussi: utilisons Proposition 2.7. \\
Surjectivit\'e: Soit $\cal E$ coh\'erent et localement libre sur $X\setminus Y$ avec une connexion et $j:X\setminus Y\to X$ l'inclusion. Alors $j_*{\cal E}$ est coh\'erent d'apr\`es Proposition 2.2. La connexion s'\'etend \`a une connexion unique sur $j_*{\cal E}$ d'apr\`es Proposition 2.7. Par cons\'equent, $j_*{\cal E}$ est localement libre.\\
b) D'apr\`es Proposition 8.1 on a $Vect(X\setminus Y)\simeq Vect(X^{an}\setminus Y^{an})$. De plus, si $\cal E$ est un fibr\'e vectoriel sur $X\setminus Y$, une connexion sur ${\cal E}^{an}$ provient d'une connexion unique sur $\cal E$: on utilise Proposition 6.2 ici.\\

Mauntenant retournons \`a la situation habituelle.

\begin{theoreme} Supposons que $X\setminus Y$ est lisse de dimension $\ge 3$, $H$ coupe $X\setminus Y$ transversalement, ${\rm codim}_XY\cap H\ge 4$. On a un diagramme commutatif:\\
$$\begin{array}{ccccc}
Vect_c(X\setminus Y)&&\simeq&&Vect_c(\hat{X}\setminus \hat{Y})\\
\downarrow\simeq&&&&\downarrow\simeq\\
Vect_c(X^{an}\setminus Y^{an})&\simeq&\lim\limits_\rightarrow Vect_c(V)&\simeq&Vect_c(\hat{X}^{an}\setminus \hat{Y}^{an})
\end{array}$$
o\`u $V$ parcourt les voisinages ouverts de $X^{an}\cap H^{an}\setminus Y^{an}$ dans $X^{an}\setminus Y^{an}$.
\end{theoreme}

{\bf D\'emonstration:} La verticale \`a gauche est bijective d'apr\`es Th\'eor\^eme 9.2b). \\
La verticale \`a droite est bijective aussi: D'abord, on a $Vect(\hat{X}\setminus\hat{Y})\simeq Vect(\hat{X}^{an}\setminus\hat{Y}^{an})$ d'apr\`es Proposition 8.2.
Soit $\cal E$ un fibr\'e vectoriel sur $\hat{X}\setminus\hat{Y}$. Supposons que nous avons une connexion sur ${\cal E}^{an}$: celle-ci peut \^etre donn\'ee par une section de $Hom(P^1({\cal E}^{an}),\Omega^1_{\hat{X}^{an}\setminus \hat{Y}^{an}}\otimes{\cal E}^{an})\simeq {\cal T}^{an}$, o\`u ${\cal T}:=Hom(P^1({\cal E}),\Omega^1_{\hat{X}\setminus \hat{Y}}\otimes{\cal E})$. Avec Proposition 6.7 on peut conclure que la section vient d'une section unique de $\cal T$, celle-ci correspond \`a une connection sur $\cal E$.\\
La horizontale sup\'erieure est injective d'apr\`es Th\'eor\^eme 4.2 et 3.1.\\
La composition des fl\`eches inf\'erieures est donc injective aussi. Surjectivit\'e: Un \'el\'ement de $Vect_c(\hat{X}^{an}\setminus\hat{Y}^{an}))$ est repr\'esent\'e par un fibr\'e ${\cal E}^{an}$ avec connexion, o\`u ${\cal E}^{an}$ provient d'un fibr\'e $\cal E$ sur $\hat{X}\setminus \hat{Y}$. D'apr\`es Th\'eor\`eme 3.6 il y a une extension coh\'erente $\cal S$ de $\cal E$ \`a $X\setminus Y$. En fait, $\dim\,\overline{S_l({\cal O}_{\hat{X}\setminus\hat{Y}}/{\cal IO}_{\hat{X}\setminus\hat{Y}})}\cap Y<l-2$ pour $l\le\dim\,Y\cap H+2$ parce que ${\rm codim}_XY\cap H\ge 4$, et $prof\,{\cal O}_{\hat{X}\setminus\hat{Y}}/{\cal IO}_{\hat{X}\setminus\hat{Y}}\ge 2$ parce que $\dim\,X\setminus Y\ge 3$. La connexion sur ${\cal E}^{an}$ est donn\'ee par une section de  $Hom(P^1({\cal E}^{an}),\Omega^1_{\hat{X}^{an}\setminus \hat{Y}^{an}}\otimes{\cal E}^{an})$. Or, la condition la condition de Proposition 6.8 est v\'erifi\'ee pour ${\cal T}_1:=Hom(P^1({\cal S}^{an}),\Omega^1_{X^{an}\setminus Y^{an}}\otimes{\cal S}^{an})$ au lieu de $\cal S$ et $m=2$: on applique [Si3] Lemma 3.15 ici. La section provient donc d'une section unique de ${\cal T}_1$, celle-ci correspond \`a une connexion sur ${\cal S}^{an}$. On conclut que ${\cal S}^{an}$ est localement libre.\\
Il reste \`a montrer que la premi\`ere fl\`eche inf\'erieure est surjective aussi. Or, on a $\lim\limits_\rightarrow\,Vect(V)\simeq Vect(X^{an}\setminus Y^{an}|X^{an}\cap H^{an}\setminus Y^{an})$ d'apr\`es Corollaire 10.4 (voir les notations y utilis\'ees). Par [Go] II Th. 3.3.1 on conclut que $\lim\limits_\rightarrow\,Vect_c(V)\simeq Vect_c(X^{an}\setminus Y^{an}|X^{an}\cap H^{an}\setminus Y^{an})$. Un \'el\'ement de $Vect_c(X^{an}\setminus Y^{an}|X^{an}\cap H^{an}\setminus Y^{an})$ est repr\'esent\'e par un $({\cal O}_{X^{an}}|X^{an}\cap H^{an}\setminus Y^{an})$-module coh\'erent localement libre $\cal E$ avec une connexion. Or, $\cal E$ peut \^etre \'etendu \`a un faisceau $\cal S$ sur $X$ d'apr\`es Th\'eor\`eme 7.4. La connexion du fibr\'e vectoriel $\cal E$ donne une section de $Hom(P^1({\cal S}),\Omega^1_X\otimes{\cal S})$ sur $X^{an}\cap H^{an}\setminus Y^{an}$. Par Th\'eor\`eme 7.1 il y a une extension unique sur $X^{an}$. La restriction \`a $X^{an}\setminus Y^{an}$ d\'efinit une connexion sur ${\cal S}|X^{an}\setminus Y^{an}$, il s'agit donc d'un fibr\'e vectoriel analytique \`a connexion.\\

En ce qui concerne la fl\`eche inf\'erieure \`a gauche, il y a une variante:

\begin{theoreme} Supposons que $X\setminus Y$ est lisse de dimension $\ge 3$. Alors il y a un isomorphisme
$$Vect_c(X^{an}\setminus Y^{an})\simeq\lim\limits_\rightarrow Vect_c(U\setminus Y^{an})$$
o\`u $U$ parcourt les voisinages ouverts de $X^{an}\cap H^{an}$ dans $X^{an}$.
\end{theoreme}

{\bf D\'emonstration:} Surjectivit\'e: Soit $U$ comme dans [Ha1] et ${\cal E}$ un fibr\'e vectoriel \`a connexion sur $U\setminus Y^{an}$. En particulier $\cal E$ est un faisceau analytique coh\'erent. D'apr\`es [Ha1] Theorem 2.2 il y a une extension coh\'erente $\cal S$ sur $X^{an}\setminus Y^{an}$. La connexion sur $\cal E$ donne une section de $Hom(P^1({\cal S},\Omega^1_{X^{an}}\otimes{\cal S}))$ sur $U\setminus Y^{an}$. Par [Ha1] Theorem 2.1 il y a une extension unique sur $X^{an}\setminus Y^{an}$ qui donne une connexion \`a ${\cal S}$. Ce faisceau est donc un fibr\'e vectoriel analytique \`a connexion.\\
Injectivit\'e: avec [Ha1] Theorem 2.1.\\

{\bf b) Fibr\'es \`a connexion int\'egrable}\\

D'abord

\begin{theoreme} Sous les hypoth\`eses de Th\'eor\`eme 9.1, on a $Vect_{ci}X\simeq Vect_{ci}(X\setminus Y)$.
\end{theoreme}

{\bf D\'emonstration:} On applique Lemme 10.9 ou bien Lemme 10.10.

\begin{theoreme} Sous les hypoth\`eses de Th\'eor\`eme 9,2, on a:\\
a) Si $X$ est lisse: $Vect_{ci}X\simeq Vect_{ci}(X\setminus Y)$, et $Vect_{cir}X\simeq Vect_{ci}(X^{an})$\\
b) Si $X$ est complet et $X\setminus Y$ lisse on a $Vect_{cir}(X\setminus Y)\simeq Vect_{ci}(X\setminus Y)\simeq Vect_{ci}(X^{an}\setminus Y^{an})$.
\end{theoreme}

{\bf D\'emonstration:} a) La premi\`ere isomorphie s'obtient comme dans la d\'emonstration de Th\'eor\`eme 9.2a). Pour la deuxi\`eme, voir [D] II Th\'eor\`eme 5.9.\\
b) Traitons d'abord l'analogue de Th\'eor\`eme 9.1. La paire $(X^{an},X^{an}\setminus Y^{an})$ est $2$-connexe, donc $\pi_k(X^{an}\setminus Y^{an},x)\simeq \pi_k(X^{an},x)$ pour $k=0,1$, $x\in X^{an}\setminus Y^{an}$.\\
Consid\'erons maintenant l'analogue de Th\'eor\`eme 9.2:\\
a) On proc\`ede comme dans la d\'emonstration de Th\'eor\`eme 9.2a).\\
b) $Vect_{cir}(X\setminus Y)\simeq Vect_{ci}(X^{an}\setminus Y^{an})$: comme en a) avec [D] loc.cit.\\
$Vect_{cir}(X\setminus Y)\simeq Vect_{ci}(X\setminus Y)$: avec la d\'edinition de la r\'egularit\'e: [D] II 4.1(i), II 4.2, parce que $codim_XY\ge 2$.

\begin{theoreme}
Soit $X\subset\mathbb{P}_N(\mathbb{C})$ un sous-espace analytique ferm\'e, $Y$ une partie analytique ferm\'ee, $X\setminus Y$ lisse, partout de dimension 
$\ge 3$, $H$ un hyperplan. Fixons une stratification de Whitney de $(X,Y)$. Supposons qu'il y a un sous-ensemble analytique $\Sigma$ de $Y\cap H$ tel que ${\rm codim}_X\Sigma\ge 3$ et $H$ coupe toutes les strates de $X$ transversalement hors de $\Sigma$. Alors on a un diagramme commutatif d'isomorphies:
$$\begin{array}{ccc}
Vect_{cir}(X\setminus Y)&\simeq&Vect_{cir}(X\cap H\setminus Y)\\
\downarrow\simeq&&\downarrow\simeq\\
Vect_{ci}(X^{an}\setminus Y^{an})&\simeq&Vect_{ci}(X^{an}\cap H^{an}\setminus Y^{an})
\end{array}$$
\end{theoreme}

{\bf D\'emonstration:} Soit $\cal E$ un fibr\'e vectoriel (analytique) sur $X^{an}\cap H^{an}\setminus Y^{an}$ \`a connexion int\'egrable. Si $X^{an}\cap H^{an}\setminus Y^{an}$ est connexe, il est donn\'e \`a isomorphie pr\`es par un homomorphisme $\pi_1(X^{an}\cap H^{an}\setminus Y^{an})\to GL_k(\mathbb{C})$ o\`u $k$ est le rang de $\cal E$. Autrement il faut regarder chaque composante connexe s\'epar\'ement. Voir [D] I.1.4, 2.17.\\
Par [HL1] on sait, si $\Sigma=\emptyset$, que pour chaque point $x$ de $X^{an}\cap H^{an}\setminus Y^{an}$, on a $\pi_k(X^{an}\cap H^{an}\setminus Y^{an},x)\simeq \pi_k(X^{an}\setminus Y^{an},x)$, $k=0,1$. En g\'en\'eral, on obtient le m\^eme \'enonc\'e, en comparant les deux vari\'et\'es d'abord avec l'intersection par un sous-espace lin\'eaire de dimension $2$ qui est contenu dans $H$ et coupe $X\setminus Y$ transversalement.\\
Ceci signifie que $\cal E$ peut \^etre \'etendu \`a un fibr\'e ${\cal F}$ sur $X^{an}\setminus Y^{an}$ \`a connexion int\'egrable qui est unique, \`a isomorphisme pr\`es. La fl\`eche inf\'erieure est donc bijective.\\
Par [D] II.5.9, ${\cal E}$ provient d'un fibr\'e vectoriel alg\'ebrique \`a connexion int\'egrable r\'eguli\`ere qui est unique, \`a isomorphie pr\`es. Le m\^eme vaut pour ${\cal F}$. Comme la restriction d'un fibr\'e vectoriel alg\'ebrique \`a connexion int\'egrable r\'eguli\`ere \`a un sous-espace lisse donne un objet de la m\^eme cat\'egorie on obtient que les fl\`eches verticales sont bijectives.

\vspace{2mm}
{\bf D\'emonstration de Th\'eor\`eme 1.8:}
Ceci r\'esulte aussit\^ot du th\'eor\`eme pr\'ec\'edant, avec $\Sigma:=Y\cap H$. Notons que les isomorphismes verticaux sont d\'ej\`a connus, par Th\'eor\`eme 9.1, on a donc $Vect_{cir}(X\setminus Y)\simeq Vect_{ci}(X\setminus Y)$ et l'analogue pour les vari\'et\'es compl\'et\'ees.\\

De nouveau il y a une variante:

\begin{theoreme}
Soit $\dim\,X\ge 3$, $X\setminus Y$ lisse. Alors il y a un diagramme commutatif
$$\begin{array}{ccc}
Vect_{cir}(X\setminus Y)&&\\
\downarrow\simeq&\searrow\simeq&\\
Vect_{ci}(X^{an}\setminus Y^{an})&\stackrel{\simeq}{\longrightarrow}&\lim\limits_\rightarrow Vect_{ci}(U\setminus Y^{an})
\end{array}$$
o\`u $U$ parcourt les voisinages ouverts de $X^{an}\cap H^{an}$ dans $X^{an}$.
\end{theoreme}

{\bf D\'emonstration:} Par [HL1] on sait que pour chaque point $x$ de $U\setminus Y^{an}$, $\pi_k(U\setminus Y^{an},x)\simeq \pi_k(X^{an}\setminus Y^{an},x)$, $k=0,1$, si $U$ est un voisinage convenable. Ceci implique la bijectivit\'e de la fl\`eche inf\'erieure.\\
La bijectivit\'e de la fl\`eche verticale se d\'emontre comme dans la d\'emonstration pr\'ec\'edante.

\section{\bf Appendice}

{\bf a) D\'emonstrations alternatives}\\

{\bf D\'emonstration directe de Th\'eor\`eme 1.2 pour $m=2$:}\\
a) D'apr\`es Proposition 2.2 on a: $R^sj_*{\cal S}$ est coh\'erent pour $s=0,1$, o\`u $j:X\setminus Y\longrightarrow X$ est l'inclusion.\\
b) D'apr\`es Lemme 2.8, le faisceau ${\cal T}:=j_*{\cal S}$ est coh\'erent, et $prof\,{\cal T}\ge 2$.\\
c) Par b) et Th\'eor\`eme 1.1 on obtient $H^s(X,j_*{\cal S}\otimes{\cal L}^{-\nu})=0$, $s=0,1, \nu\gg 0$.\\
d) Soit $x$ un point ferm\'e de $X$. D'apr\`es b) on a $H^1_x(j_*{\cal S})=0$, donc $\mathbb{H}^1_x(Rj_*{\cal S})\simeq H^0_x(R^1j_*{\cal S})$.\\
Mais $\mathbb{H}^1_x(Rj_*{\cal S})\simeq \lim\limits_\rightarrow \mathbb{H}^1(U,U\setminus\{x\},Rj_*{\cal S})\simeq \lim\limits_\rightarrow \mathbb{H}^1(U\setminus Y,U\setminus(Y\cup \{x\}),{\cal S})=0$ si $x\in Y$. En plus, $R^1j_*{\cal S}$ est concentr\'e sur $Y$. En total, $H^0_x(R^1j_*{\cal S})=0$ toujours, donc $prof\,R^1j_*{\cal S}\ge 1$. \\
e) Par d) et Th\'eor\`eme 1.1 on conclut que $H^0(X,(R^1j_*{\cal S})\otimes{\cal L}^{-\nu})=0$, $\nu\gg 0$.\\
f) Par la formule d'adjonction on obtient, ${\cal L}^{-\nu}$ \'etant localement libre: \\
$(R^sj_*{\cal S})\otimes{\cal L}^{-\nu}\simeq R^sj_*({\cal S}\otimes{\cal L}^{-\nu})$. \\
Par c), e) on obtient donc que
$H^p(X,R^qj_*({\cal S}\otimes{\cal L}^{-\nu}))=0$, $p+q\le 1,\nu\gg 0$. Ceci implique notre \'enonc\'e.\\

{\bf Autre d\'emonstration de Lemme 2.5:}\\
On peut supposer que $\cal S$ est coh\'erent sur $X$. Alors
$\oplus_{n\ge 0} H^s(X,{\cal S}(n))$ est un $k[Z_0,\ldots,Z_r]$-module de type fini pour tout $s$: Pour $s=0$ voir [G4] III 2.3.2, pour $s>0$ on
sait que $H^s(X,{\cal S}(n))$ est de dimension fini et $0$ si $n\gg 0$.\\
On en d\'eduit que $\oplus_{n\ge 0} H^s(X\setminus Y,{\cal S}(n))$ est un $k[Z_0,\ldots,Z_r]$-module de type fini pour $s<m$: il suffit de montrer que $\oplus_{n\ge 0} H^s(X,R^lj_*{\cal S}(n))$ est un $k[Z_0,\ldots,Z_r]$-module de type fini pour $l<m$, et $R^lj_*{\cal S}$ est coh\'erent d'apr\`es Proposition 2.2. \\
D'autre part, $H^s(X\setminus Y,{\cal S}(n))$ est de dimension finie pour $n$ fixe et $s<m$, par Proposition 2.2, et $H^s(X\setminus Y,{\cal S}(n))=0$
pour $n\ll 0$ d'apr\`es Th\'eor\`eme 1.2. D'o\`u le r\'esultat.\\

{\bf D\'emonstration directe de Corollaire 3.7:}\\
Soit $\hat{j}:\hat{X}\setminus \hat{Y}\to \hat{X}$ l'inclusion. D'apr\`es Proposition 2.2, les faisceaux $j_*({\cal S}/{\cal IS})$ et $R^1j_*({\cal S}/{\cal IS})$ sont coh\'erents. Par cons\'equent, $\oplus_{k\ge 0}R^lj_*({\cal I}^k{\cal S}/{\cal I}^{k+1}{\cal S})$ est un $\oplus_{k\ge 0}{\cal I}^k/{\cal I}^{k+1}$-module de type fini, $l=0,1$. D'apr\`es [G2] IX Th. 2.1 on peut conclure que ${\cal T}:=\hat{j}_*{\cal S}$ est coh\'erent. ce qui implique Corollaire 3.7.\\

{\bf Remarque:} La suite exacte\\
$0\to {\cal H}^0_{\hat{Y}}{\cal T}\to {\cal T}\stackrel{\simeq}{\to} \hat{j}_*\hat{j}^*{\cal T}\to {\cal H}^1_{\hat{Y}}{\cal T}\to 0$\\
implique que ${\cal H}^l_{\hat{Y}}{\cal T}=0, l\le 1$, donc $prof\,{\cal T}_y\ge 2$ pour $y\in\hat{Y}$. On a donc $prof\,{\cal T}\ge 2$. Mais ceci ne suffit pas afin d'appliquer [G2] XII Th\'eor\`eme 3.1 pour une extension coh\'erente de $\cal T$ \`a $X$ et d\'emontrer Th\'eor\`eme 3.6!\\
On devrait m\^eme avoir $prof\,{\cal T}/{\cal IT}\ge 2$, c.-\`a d. ${\cal T}/{\cal IT}\simeq j_*j^*({\cal T}/{\cal IT})$.\\

{\bf D\'emonstration directe de Th\'eor\`eme 3.6:}\\
Notons que $\cal S$ induit un faisceau $\tilde{\cal S}$ sur $\hat{\tilde{X}}\setminus\hat{\tilde{Y}}$. Soit $\hat{\tilde{j}}:\hat{\tilde{X}}\setminus\hat{\tilde{Y}}\to \hat{\tilde{X}}$ l'inclusion. D'apr\`es [G2] IX Th. 2.1, ${\cal T}:=\hat{\tilde{j}}_*\hat{\tilde{\cal S}}$ est coh\'erent.\\
On continue maintenant comme dans la d\'emonstration de la surjectivit\'e de $\lim\limits_\rightarrow Vect\,U\to Vect\,\hat{X}\setminus\hat{Z}$ \`a la fin de la d\'emonstration de Th\'eor\`eme 7.2 dans [HL2].\\

Ceci devrait revenir essentiellement \`a une d\'emonstration du r\'esultat de Gro\-thendieck [G2] XII Th\'eor\`eme 3.1 dans le sens de sa d\'emonstration originale comme il l'y indique (p. 149).\\

{\bf D\'emonstration alternative de Lemme 4.4:}\\
Soit $k\ge n$ et $\cal E$ un fibr\'e vectoriel sur $X_{k+1}$ pour lequel la restriction \`a $X_k$ est trivial. Soient $s_1,\ldots,s_r\in \Gamma(X\cap H\setminus Y,{\cal E}|X_k)$ des sections qui donnent une base pour chaque fibre. L'hypoth\'ese implique que les sections s'\'etendent \`a $X_{k+1}$, elles d\'efinissent la m\^eme base pour chaque fibre; en fait, on a une suite exacte\\
$\Gamma(X\cap H\setminus Y,{\cal E}/{\cal I}^{k+1}{\cal E})\to\Gamma(X\cap H\setminus Y,{\cal E}/{\cal I}^k{\cal E})\to H^1(X\cap H\setminus Y,{\cal I}^k{\cal E}/{\cal I}^{k+1}{\cal E})$\\
et ${\cal I}^k{\cal E}/{\cal I}^{k+1}{\cal E}\simeq ({\cal I}^k/{\cal I}^{k+1})\otimes({\cal E}/{\cal I}{\cal E})\simeq ({\cal I}^k/{\cal I}^{k+1})^r$, o\`u $r$ est le rang de $\cal E$, parce que ${\cal E}/{\cal IE}$ est trivial. On a donc $H^1(X\cap H\setminus Y,{\cal I}^k{\cal E}/{\cal I}^{k+1}{\cal E})=0$ d'apr\`es l'hypoth\`ese.\\

{\bf Autre d\'emonstration de Th\'eor\`eme 1.7:}\\
Autre m\'ethode: Supposons que ${\cal E}/{\cal IE}$ est trivial. Soient $s_1,\ldots,s_k$ des section de ce fibr\'e qui donnent une base pour chaque fibre. Maintenant $({\cal I}^n{\cal E})/({\cal I}^{n+1}{\cal E})\simeq {\cal I}^n/{\cal I}^{n+1})\otimes({\cal E}/{\cal IE})\simeq ({\cal I}^n/{\cal I}^{n+1}$, donc $H^l(X\cap H\setminus Y,({\cal I}^n{\cal E})/({\cal I}^{n+1}{\cal E}))=0, l=0,1$ \`a cause de Th\'eor\`eme 8.5. Ceci implique que $s_1,\ldots,s_k$ proviennent des sections de ${\cal E}/{\cal I}^n{\cal E}$, $n\gg 0$, donc des sections de $\cal E$, par Th\'eor\`eme 1.2. Celles-ci d\'efinissent encore une base de chaque fibre sur $X\cap H\setminus Y$, donc ${\cal E}$ est trivial sur un voisinage ouvert de Zariski de $X\cap H\setminus Y$ dans $X\setminus Y$. Par Th\'eor\`eme 4.1 $\cal E$ est trivial.\\

{\bf Deuxi\`eme d\'emonstration de Th\'eor\`eme 1.8:}\\
On utilise d'abord Th\'eor\`eme 9.1 o\`u on peut remplacer $c$ par $ci$ partout. Car l'int\'egrabilit\'e d'une connexion sur $\cal E$ signifie que la courbure est $0$, la courbure \'etant une section de $Hom({\cal E},\Omega^2_X\otimes_{{\cal O}_X}{\cal E})$. On peut donc raisonner avec Proposition 6.8.\\
On a aussi $Vect_{ci}(X\cap H\setminus Y)\simeq Vect_{ci}(X^{an}\cap H^{an}\setminus Y^{an})$. D'autre part, la d\'emonstration pr\'ec\'edante montre qu'on peut \'etendre chaque fibr\'e \`a connexion int\'egrable sur $X\cap H\setminus Y$ \`a un voisinage, donc \`a $\hat{X}\setminus\hat{Y}$. Par cons\'equent, $Vect_{ci}(X\setminus Y)\to Vect_{ci}(X\cap H\setminus Y)$ est surjectif.\\
Injectivit\'e: Soient ${\cal E}_1, {\cal E}_2$ des fibr\'es vectoriels sur $X^{an}\setminus Y^{an}$ \`a connexion int\'egrable dont la restriction \`a $X^{an}\cap H^{an}\setminus Y^{an}$ est isomorphe. Alors la restriction \`a un voisinage est isomorphe aussi, donc $\hat{\cal E}_1\simeq\hat{\cal E}_2$. Par Th\'eor\`eme 8.7 on obtient l'isomorphie ${\cal E}_1\simeq{\cal E}_2$.\\ 

{\bf D\'emonstration alternative de Th\'eor\`eme 9.5:} On proc\`ede comme dans la d\'emonstration de Th\'eor\`eme 9.1.\\

{\bf D\'emonstration alternative de Th\'eor\`eme 9.6 b):} On d\'emontre $Vect_{ci}(X\setminus Y)\simeq Vect_{ci}(X^{an}\setminus Y^{an})$ par voie directe, comme Th\'eor\`eme 9.2 b).\\

{\bf b) Compl\'ements sur la th\'eorie des faisceaux}\\

Dans ce paragraphe soit $X$ un espace paracompact, $A$ une partie ferm\'ee. Soit $(U_i)_{i\in I}$ un recouvrement ouvert localement fini de $A$, $(V_i)_{i\in I}$ un recouvrement pareil tel que $\overline{V_i}\subset U_i, i\in I$, ${\cal T}_i$ un faisceau sur $U_i$, $i\in I$. Ici on travaille avec des faisceaux d'ensembles. 

\begin{theoreme} Soit $\cal S$ un faisceau sur $A$, ${\cal S}|A\cap U_i\simeq {\cal T}_i|A\cap U_i$, $i\in I$, alors il y a un voisinage ouvert $U$ de $A$ dans $X$ et un faisceau $\tilde{\cal S}$ sur $U$ tel que $\tilde{\cal S}|U\cap V_i\simeq{\cal T}_i|U\cap V_i, i\in I$.
\end{theoreme}

{\bf D\'emonstration:} On proc\`ede comme dans la d\'emonstration de [Go] II Th. 3.3.1. Fixons des isomorphismes ${\cal S}|A\cap U_i\simeq {\cal T}_i|A\cap U_i$. Ceci donne un cocycle $(g_{ij})$ par rapport au recouvrement $(A\cap \overline{V_i})$ de $A$, o\`u $g_{ij}\in 
\Gamma(A\cap \overline{V_i}\cap \overline{V_j},Iso({\cal T}_j|U_i\cap U_j,{\cal T}_i|U_i\cap U_j))$. Par [Go] II Th. 3.3.1 il y a une extension de la section $g_{ij}$ \`a une section $\tilde{g}_{ij}$ qui est d\'efinie sur un voisinage ouvert $W_{ij}$ de $A\cap \overline{V_i}\cap \overline{V_j}$. Soit $U$ l'ensemble de tous les $x\in X$ avec les propri\'et\'es suivantes:\\
si $x\in \overline{V_i}\cap \overline{V_j}$, on a $x\in W_{ij}$,\\
si $x\in \overline{V_i}\cap \overline{V_j}\cap\overline{V_k}$, on a $\tilde{g}_{ij}(x)\circ \tilde{g}_{jk}(x)=\tilde{g}_{ik}(x)$.\\
Alors $U$ est ouvert: Soit $x\in U$. Il y a un voisinage $W(x)$ de $x$ qui ne rencontre qu'un nombre fini des $\overline{V}_i$, disons $\overline{V}_{i_1},\ldots,\overline{V}_{i_n}$. En plus, on peut supposer que $x\in \overline{V_i}$ pour $i=i_1,\ldots,i_n$. Maintenant on peut achever que $W(x)\subset U$.\\
\'Evidemment, $U$ est un voisinage de $A$, il s'agit du voisinage cherch\'e: on recolle les ${\cal T}_i|U\cap V_i$ par $\tilde{g}_{ij}|V_i\cap V_j\cap U$.

\begin{corollaire} Soit $A$ une partie analytique ferm\'ee de l'espace analytique complexe $X$. Chaque ${\cal O}_X|A$-module coh\'erent (resp. coh\'erent localement libre) admet une extension \`a un ${\cal O}_U$-module coh\'erent (resp. coh\'erent localement libre), $U$ \'etant un voisinage ouvert convenable de $A$ dans $X$.
\end{corollaire}

Pour l'unicit\'e, on a

\begin{theoreme} Si ${\cal S}_1$ et ${\cal S}_2$ sont des faisceaux sur $X$ tels que ${\cal S}_1|A\simeq {\cal S}_2|A$ alors il y a un voisinage $U$ de $A$ tel que ${\cal S}_1|U\simeq {\cal S}_2|U$.
\end{theoreme}

{\bf D\'emonstration:} On a une section du faisceau $Iso({\cal S}_1,{\cal S}_2)$ sur $A$, elle peut \^etre \'etendue sur un voisinage.

\begin{corollaire} Soit $A$ une partie analytique ferm\'ee de l'espace analytique complexe $X$, alors 
$$\lim\limits_\rightarrow Coh\,U\simeq Coh(X|A)\,,\, \lim\limits_\rightarrow Vect\,U\simeq Vect(X|A)$$
Ici $U$ parcourt tous les voisinages ouverts de $A$ dans $X$, et $Coh\,U$ resp. $Coh(X|A)$ d\'enote la cat\'egorie des ${\cal O}_U$- resp. ${\cal O}_X|A$-modules coh\'erents, \`a isomorphie pr\`es. Finalement, $Vect(X|A)$ d\'enote la cat\'egorie des ${\cal O}_X|A$-modules coh\'erents localement libres, \`a isomorphie pr\`es.
\end{corollaire}

La deuxi\`eme partie du corollaire peut aussi \^etre \'ecrite comme
$$\lim\limits_\rightarrow H^1(U,GL({\cal O}_X))\simeq H^1(A,GL({\cal O}_X))$$
Cette relation est en m\^eme temps cons\'equence d'un th\'eor\`eme g\'en\'eral sur la cohomologie des faisceaux non-ab\'eliens; dans le cas ab\'elien voir [Go] II Th. 4.11.1:

\begin{theoreme} Soit ${\cal F}$ un faisceau de groupes (non-n\'ecessairement ab\'eliens) sur $X$. Alors
$$\lim\limits_\rightarrow H^1(U,{\cal F})\simeq H^1(A,{\cal F})$$
o\`u $U$ parcourt tous les voisinages ouverts de $A$ dans $X$.
\end{theoreme}

{\bf D\'emonstration:} Surjectivit\'e: Supposons qu'un \'el\'ement de $H^1(A,{\cal F})$ est donn\'e. Nous povons supposer qu'il est repr\'esent\'e par un cocycle $(g_{ij})$ par rapport au recouvrement 
$(A\cap \overline{V_i})$ de $A$. On peut \'etendre $g_{ij}$ \`a un voisinage $W_{ij}$ de $A\cap\overline{V_i}\cap\overline{V_j}$. On continue comme dans la d\'emonstration de Th\'eor\`eme 10.1.\\
Injectivit\'e: Soient $(g_{ij})$ et $(h_{ij})$ des cocycles par rapport \`a $(U_i)$ qui d\'efinissent le m\^eme \'el\'ement de $H^1(A,{\cal F})$. Quitte \`a r\'etr\'ecir il y a $f_i\in \Gamma(A\cap\overline{V_i},{\cal F})$ tels que $h_{ij}=f_ig_{ij}f_j^{-1}$ sur $A\cap\overline{V_i}\cap\overline{V_j}$. On peut \'etendre $f_i$ \`a une section $\tilde{f}_i$ sur un voisinage ouvert $W_i$ de $A\cap\overline{V_i}$ dans $U_i$.
Soit $U$ l'ensemble de tous les points $x$ avec les propri\'et\'es suivantes:\\
si $x\in \overline{V_i}$ alors $x\in W_i$,\\
si $x\in \overline{V_i}\cap \overline{V_j}$ alors $h_{ij}(x)=\tilde{f}_i(x)g_{ij}(x)\tilde{f}_j(x)^{-1}$.\\
On v\'erifie que $U$ est un voisinage ouvert de $A$ dans $X$ et que $(h_{ij}|U\cap V_i\cap V_j)$ et $(g_{ij}|U\cap V_i\cap V_j)$ sont cohomologues, il d\'efinissent donc le m\^eme \'el\'ement de $H^1(U,{\cal F})$.\\

{\bf c) Connexions}\\

Soit $X$ une vari\'et\'e amalytique complexe, $\cal E$ un fibr\'e vectoriel holomorphe sur $X$. Il y a plusieurs fa\c{c}ons de d\'efinir la notion d'une connexion sur $X$, voir [D], [Ma]:

\vspace{2mm}
Soit ${\cal J}$ l'id\'eal de la diagonale dans $X^2$, $p_1,p_2=$ projections canoniques de $X^2$ sur $X$, $P^1:={\cal O}_{X^2}/{\cal J}^2$ et $P^1({\cal E}):=(p_1)_*(P^1\otimes_{{\cal O}_{X^2}}p_2^*{\cal E})$ le faisceaux des $1$-jets dans $\cal E$. La suite exacte\\
$0\to {\cal J}/{\cal J}^2\to {\cal O}_{X^2}/{\cal J}^2\to {\cal O}_{X^2}/{\cal J}\to 0$\\
conduit \`a une suite exacte
$$0\rightarrow \Omega^1_X\otimes{\cal E}\stackrel{i}{\rightarrow}P^1({\cal E})\stackrel{q}{\rightarrow}{\cal E}\rightarrow 0$$
Notons que $(p_1)_*(({\cal J}/{\cal J}^2)\otimes_{{\cal O}_{X^2}}p_2^*{\cal E})\simeq (p_2)_*(({\cal J}/{\cal J}^2)\otimes_{{\cal O}_{X^2}}p_2^*{\cal E})\simeq \Omega^1_X\otimes{\cal E}$ parce que $(p_2)_*({\cal J}/{\cal J}^2)\simeq \Omega^1_X$, et $(p_1)_*(({\cal O}_{X^2}/{\cal J})\otimes_{{\cal O}_{X^2}}p_2^*{\cal E}\simeq (p_2)_*(({\cal O}_{X^2}/{\cal J})\otimes_{{\cal O}_{X^2}}p_2^*{\cal E})\simeq {\cal E}$ parce que $(p_2)_*({\cal O}_{X^2}/{\cal J})\simeq {\cal O}_X$.\\
On a une application $j^1:{\cal E}\rightarrow P^1({\cal E})$ qui associe \`a chaque section son $1$-jet, induite par $f\mapsto f\circ p_2$. Elle est seulement $\mathbb{C}_X$-lin\'eaire: on a $j^1(fs)=i(df\otimes s)+fj^1(s)$ pour $f\in{\cal O}_{X,x}$, $s\in{\cal E}_x$. On a $q\circ j^1=id$, donc $q$ est vraiment surjectif, et $j^1$ donne alors un scindage de la suite exacte dans la cat\'egorie des $\mathbb{C}_X$-modules.

\begin{proposition}
Il est \'equivalent de donner\\ 
a) une application $\mathbb{C}_X$-lin\'eaire $\nabla:{\cal E}\to\Omega^1_X\otimes{\cal E}$ telle que $\nabla(fs)=df\otimes s+f \nabla s$,\\
b) une application ${\cal O}_X$-lin\'eaire $D_1:P^1({\cal E})\to \Omega^1_X\otimes{\cal E}$ telle que $D_1\circ i=id$,\\
c) une application ${\cal O}_X$-lin\'eaire $D_2:{\cal E}\to P^1({\cal E})$ telle que $q\circ D_2=id$,\\
o\`u les applications sont li\'ees comme suit: $D_2:=j^1-i\circ\nabla$, $i\circ D_1=id-D_2\circ q$, $\nabla:=D_1\circ j^1$.
\end{proposition}

Pour $a\Leftrightarrow c)$ voir [D), o\`u on trouve la relation $D_2=j^1-i\circ\nabla$ aussi. Notons que $c)$ signifie que $D_2$ donne un scindage de la suite exacte en haut dans la cat\'egorie des ${\cal O}_X$-modules. Il est bien connu que ceci est \'equivalent \`a $b)$ avec la relation $D_2\circ q+i\circ D_1=id$. On obtient $\nabla=D_1\circ j^1$ comme cons\'equence; cette relation semble \^etre conforme avec [Ma]. En total, on voit comme on passe de $\nabla$ \`a $D_2$, de $D_2$ \`a $D_1$ et de $D_1$ \`a $\nabla$. Pour compl\'eter les relations signalons qu'on peut raisonner au sens inverse aussi: $i\circ D_1=i\circ\nabla\circ q+id-j^1\circ q$, $D_2=j^1-i\circ D_1\circ j^1$, $i\circ\nabla=j^1-D_2$.

\vspace{2mm}
La premi\`ere d\'efinition est la d\'efinition usuelle mais les autres ont un avantage du point de vue th\'eorique parce qu'il s'agit des applications ${\cal O}_X$-lin\'eaires.\\ 

{\bf Remarque:} On peut remplacer $\cal E$ par un faisceau coh\'erent. La suite en haut reste exacte parce qu'on a encore $q\circ j^1=id$, $q$ est donc surjectif. Proposition 10.6 reste encore valable. Mais par [Ma] Remark (1.2), l'existence de $\nabla$ garantit que le faisceau coh\'erent est automatiquement localement libre.\\

Rappelons qu'une connexion $\nabla$ sur $\cal E$ induit des applications $\mathbb{C}_X$-lin\'eaires $\nabla^p:\Omega^p_X\otimes{\cal E}\to\Omega^{p+1}_X\otimes{\cal E}$, o\`u $\nabla^0=\nabla$, voir [D]. Or,  $\nabla$ est dit int\'egrable si $\nabla^1\circ\nabla=0$.\\

En ce qui concerne les connexions sur un fibr\'e vectoriel donn\'e, on a le r\'esultat suivant:

\begin{proposition}
Soit $\cal E$ un fibr\'e vectoriel sur $X$.\\
a) Si $H^0(X,\Omega^1_X\otimes Hom({\cal E},{\cal E}))=0$ il y a au plus une connexion sur $\cal E$.\\
b) Si $H^1(X,\Omega^1_X\otimes Hom({\cal E},{\cal E}))=0$ il y a au moins une connexion sur $\cal E$.\\
c) Si $H^0(X,\Omega^2_X\otimes Hom({\cal E},{\cal E}))=0$ toute connexion sur $\cal E$ est int\'egrable.
\end{proposition}

{\bf D\'emonstration:} a), b): La suite exacte (*) en haut conduit \`a une suite exacte
$$0\to Hom({\cal E},\Omega^1_X\otimes{\cal E})\to Hom({\cal E},P^1({\cal E}))\to Hom({\cal E},{\cal E})\to 0$$
donc\\
$H^0(X,Hom({\cal E},\Omega^1_X\otimes{\cal E}))\to H^0(X,Hom({\cal E},P^1({\cal E})))\to H^0(X,Hom({\cal E},{\cal E}))$\\
$\to H^1(X,Hom({\cal E},\Omega^1_X\otimes{\cal E}))$.\\
Or, $id:{\cal E}\to{\cal E}$ d\'efinit un \'el\'ement de $H^0(X,Hom({\cal E},{\cal E}))$, et une connexion sur $\cal E$ correspond \`a une image inverse dans $H^0(X,Hom({\cal E},P^1({\cal E})))$, par Proposition 10.6. Le reste est clair.\\
c) Soit $\nabla^1$ d\'efini comme toute \`a l'heure. Or, l'application $\nabla^1\circ\nabla$ est m\^eme ${\cal O}_X$-lin\'eaire, voir [D], elle d\'efinit donc un \'el\'ement de $H^0(X,\Omega^2_X\otimes Hom({\cal E},{\cal E}))$.\\

Soit $Vect_c(X)$ et $Vect_{ci}(X)$ d\'efini comme dans l'introduction. D'apr\`es [D] nous avons que $Vect_c(X)$ a la structure d'un semi-anneau avec unit\'e (donn\'ee par $({\cal O}_X,d)$), le m\^eme vaut pour $Vect_{ci}(X)$. De plus, on a $Vect_{ci}(X)\simeq H^1(X,GL(\mathbb{C}))$.\\

{\bf d) Un lemme sur les suites spectrales}\\

Le lemme suivant devrait \^etre bien connu (il est implicitement appliqu\'e dans [H2] p. 223, par exemple):

\begin{lemme}
Soit $f:E\to E'$ un homomorphisme de suites spectrales $(E_r)_{r\ge r_0}$, $(\tilde{E}_r)_{r\ge r_0}$, avec $r_0\ge 1$. Supposons que $E^{pq}_{r_0}=\tilde{E}^{pq}_{r_0}=0$ si $p<0$ ou $q<0$ et que $f$ induit des isomorphismes $E^{pq}_{r_0}\simeq \tilde{E}^{pq}_{r_0}$ si $q\le n$ (il suffit: si $q\le n-\frac{r_0-1}{r_0}p$). Alors $f$ induit des isomorphismes $E^{pq}_{\infty}\simeq\tilde{E}^{pq}_{\infty}, p+q\le n$.
\end{lemme}

{\bf D\'emonstration:} On d\'emontre par r\'ecurrence que $f$ induit des isomorphismes $E^{pq}_r\simeq\tilde{E}^{pq}_r, q\le n-\frac{r-1}{r}p$, $r\ge r_0$: afin de d\'emontrer cet \'enonc\'e avec $r+1$ au lieu de $r$ utilisons le diagramme commutatif
$$\begin{array}{ccccc}
E_r^{p-r,q+r-1}&\longrightarrow&E_r^{p,q}&\longrightarrow&E_r^{p+r,q-r+1}\\
\downarrow\simeq&&\downarrow\simeq&&\downarrow\simeq\\
\tilde{E}_r^{p-r,q+r-1}&\longrightarrow&\tilde{E}_r^{p,q}&\longrightarrow&\tilde{E}_r^{p+r,q-r+1}
\end{array}$$
pour $q\le n-\frac{r-1}{r}p$.\\
Notons aussi que $E^{pq}_r=E^{pq}_\infty$ et $\tilde{E}^{pq}_r =\tilde{E}^{pq}_\infty$ si $p+q\le n, r>n$.\\

{\bf e) Sur les groupes d'homotopie et syst\`emes locaux}\\

Soit $X$ une vari\'et\'e analytique r\'eelle et $Y$ un sous-ensemble analytique ferm\'e de $X$. Le lemme suivant est g\'eom\'etriquement \'evident, nous fournissons une preuve qui n'utilise pas des arguments de transversalit\'e:

\begin{lemme} 
Supposons que $d$ est la codimension de $Y$ dans $X$. Alors la paire $(X,X\setminus Y)$ est $d-1$-connexe, c.-\`a d. $\pi_l(X\setminus Y,x)\to\pi_l(X,x)$, est bijectif pour $l\le d-2$ et surjectif pour $l=d-1$, $x\in X\setminus Y$.
\end{lemme}

{\bf D\'emonstration:} On peut supposer que $X$ est connexe de dimension $n$. Il y a une filtration $\emptyset=Y_{-1}\subset Y_0\subset Y_1\subset\ldots\subset Y_{n-d}=Y$ o\`u $Y_k$ est analytique ferm\'e de dimension $k$ et $Y_k\setminus Y_{k-1}$ est lisse de dimension $k$. Il suffit de d\'emontrer que $\pi_l(X\setminus Y_k,x)\simeq \pi_l(X\setminus Y_{k-1},x)$, $x\in X\setminus Y_k, l\le n-k-2$. En rempla\c{c}ant $X\setminus Y_{k-1}$ par $X$ on voit qu'il suffit de d\'emontrer notre \'enonc\'e original dans le cas $Y$ lisse.\\
Soit $U$ un voisinage tubulaire de $Y$ dans $X$, on a donc une r\'etraction $r:U\to Y$. Il suffit de montrer $\pi_l(U\setminus Y,x)\simeq\pi_l(U,x)$, $l\le d-2$, $x\in U\setminus Y$, \`a cause du th\'eor\`eme de Blakers-Massey. Soit $F:=r^{-1}(\{x\})$; alors $F$ est hom\'eomorphe \`a un disque de dimension $n-d$. Or, $r|U\setminus Y\to Y$ d\'efinit une fibration; la suite exacte d'homotopie pour celle-ci donne: $\pi_l(U\setminus Y,x)\simeq \pi_l(Y,x)$, donc  $\pi_l(U\setminus Y,x)\simeq \pi_l(U,x)$ pour $l\le d-2$, ce qu'il fallait d\'emontrer.

\begin{lemme}
Soit $\cal L$ un syst\`eme local sur $X\setminus Y$ et $Y$ de codimension $\ge 3$. Alors $\cal L$ admet comme syst\`eme local une extension \`a $X$ qu est unique \`a isomorphie pr\`es, \`a savoir $j_*{\cal L}$, o\`u $j:X\setminus Y\to X$ est l'inclusion.
\end{lemme}

{\bf D\'emonstration:} Avec les notations de la d\'emonstration pr\'ec\'edente, il suffit de d\'emontrer par r\'ecurrence que $j_*{\cal L}$ est un syst\`eme local sur $X\setminus Y_k$. On peut donc supposer que $Y$ est lisse. Il suffit de montrer que $j_*{\cal L}|U$ est un syst\`eme local, $U$ \'etant un voisinage convenable de $x\in Y$ dans $X$, ce qui est \'evident.\\
Ou bien on applique Lemme 10.9.

\vspace{2cm}

{\bf Bibliographie:}\\

[BS] C. B\u anic\u a, O. St\u an\u a\c sil\u a, Algebraic methods in the 
Global Theory of Complex Spaces, John Wiley, London, N.Y., Toronto, 1976.\\

[D] P. Deligne: \'Equations diff\'erentielles \`a points  singuliers r\'eguliers. Springer Lecture Notes in Math. {\bf 163} (1970).\\

[FG] J. Frisch, J. Guenot: Prolongement de faisceaux analytiques coh\'erents. Invent. Math. {\bf 7}, 321-343 (1969).\\

[Fr] J. Frenkel: Cohomologie non ab\'elienne et espaces fibr\'es. Bull. Soc. Math. France {\bf 85}, 135-220 (1957).\\

[G1] A. Grothendieck: El\'ements de la G\'eom\'etrie I. Publ. Math. IHES {\bf 4} (1960).\\

[G2] A. Grothendieck: Cohomologie locale des faisceaux coh\'erents et 
th\'eor\`emes de Lefschetz locaux et globaux (SGA II). North-Holland: 
Amsterdam 1968.\\

[G3] A. Grothendieck, J.A.Dieudonn\'e: \'El\'ements de G\'eom\'etrie Alg\'ebrique I. Springer: New York 1971.\\

[G4] A. Grothendieck: El\'ements de la G\'eom\'etrie III, premi\`ere partie. Publ. Math. IHES {\bf 11} (1961).\\

[G5] A. Grothendieck: El\'ements de la G\'eom\'etrie III, seconde partie. Publ. Math. IHES {\bf 17} (1963).\\

[Go] R. Godement: Topologie alg\'ebrique et th\'eorie des faisceaux. Hermann: Paris 1958.\\

[GR] R. Gunning, H.Rossi: Analytic Functions of Several Complex Variables. Prentice-Hall: Englewood Cliffs, N.J. 1965.\\

[H1] R. Hartshorne: Algebraic geometry. Springer: New York 1977.\\

[H2] R. Hartshorne: Ample Subvarieties of Algebraic Varieties. Springer Lecture Notes in Math. {\bf 156} (1970).\\

[Ha1] H.A. Hamm: On theorems of Zariski-Lefschetz type. In: Singularities II: Geometrical and topological aspects, pp. 69-78, ed. J.-P. Brasselet et al. Contemporary Math. {\bf 475}. Amer. Math. Soc., Providence, R.I. 2008.\\

[Ha2] H.A. Hamm: On the local Picard group. Proc. Steklov Inst. of Math. {\bf 267}, 131-138 (2009).\\

[Ha3] H.A. Hamm: Chow groups and tubular neighbourhoods. J. Singul. {\bf 2}, 81-91 (2010).\\

[HL1] H.A. Hamm, L\^e D.T.: Lefschetz theorems on quasi-projective varieties. Bull. Soc. Math. France {\bf 113}, 123-142 (1985).\\

[HL2] H.A. Hamm, L\^e D.T.: Th\'eor\`emes d'annulation et groupes de Picard. J. Singul. {\bf 1}, 13-36 (2010).\\

[KK] L. Kaup, B. Kaup: Holomorphic functions of several variables. Walter de 
Gruyter: Berlin 1983.\\

[M] H. Matsumura: Commutative ring theory. Cambridge Univ. Press: Cambridge 1986.\\

[Ma] B. Malgrange: Regular Connections after Deligne. In: A.Borel et al.: Algebraic $D$-modules, Ch. IV, pp. 151-172. Perspectives in Math. Vol. 2. Acad. Press: Boston 1987.\\

[Sch] G. Scheja: Fortsetzungss\"atze der komplex-analytischen Cohomologie und ihre algebraische Charakterisierung. Math. Ann. {\bf 57}, 75-94 (1964).\\

[Si1] Y.-T. Siu: Extending coherent analytic sheaves. Ann. of Math. {\bf 90}, 108-143 (1969).\\

[Si2] Y.-T. Siu: Analytic sheaves of local cohomology. Trans. Am. Math. Soc. {\bf 148}, 347-366 (1970).\\

[Si3] Y.-T. Siu: Techniques of extension of analytic objects. Lecture Notes in Pure and Applied Math. {\bf 8}. Marcel Dekker: N.Y. 1974.\\

[SiT] Y.-T. Siu, G. Trautmann: Gap-sheaves and extensions of coherent analytic subsheaves. Springer Lecture Notes in Math. {\bf 172} (1971).\\

[T1] G. Trautmann: Ein Kontinuit\"atssatz f\"ur die Fortsetzung koh\"arenter analytischer Garben. Archiv der Math. {\bf 18}, 188-196 (1967).\\

[T2] G. Trautmann: Ein Endlichkeitssatz in der analytischen Geometrie. Invent. Math. {\bf 8}, 143-174 (1969).\\

\end{document}